\documentclass[12pt]{article}

\usepackage[utf8]{inputenc}
\usepackage{graphicx}
\usepackage{enumitem}
\usepackage{amsmath, bm, amsthm, esint}
\usepackage{listings}
\usepackage{hyperref}
\usepackage{float}
\usepackage{comment}
\usepackage{xcolor}
\usepackage{subcaption}
\usepackage[title]{appendix}
\usepackage{makecell}
\usepackage{authblk}
\usepackage[symbol]{footmisc}

\usepackage{lineno}

\usepackage[sorting=none, style=numeric-comp]{biblatex} 
\newtheorem*{remark}{Remark}

\usepackage[left=1.00 in, right=1.00 in, top=1.00 in, bottom=1.00 in]{geometry}
\title{\LARGE \textbf{A Path-Conservative Method for a Weakly Compressible Two-Phase Model with Surface Tension}}
\author{Ashley Melvin, J. C. Mandal\footnote{Corresponding author: \href{mailto:mandal@aero.iitb.ac.in}{mandal@aero.iitb.ac.in} (J. C. Mandal).}}
\affil{\small Department of Aerospace Engineering, Indian Institute of Technology Bombay, Mumbai, Maharashtra, 400076, India}
\date{}

\begin{document}

\maketitle

\begin{abstract}
	When extended to two-phase flows, weakly compressible models lead to a non-conservative system, which precludes its treatment using standard finite volume techniques. In this paper, a novel HLLC-type path-conservative scheme is formulated for the weakly compressible two-phase model. Furthermore, capillary effects are included in the proposed path-conservative scheme, eliminating the need to discretize surface tension terms separately. The first-order path-conservative formulation is combined with a local solution reconstruction technique to obtain high-order spatial accuracy. The solver is tested on several benchmark two-phase flow problems to demonstrate its efficacy.
\end{abstract}



\section{Introduction}\label{sec:into}

Riemann solvers have proven to be a reliable method to solve incompressible two-phase flows \cite{yang2014upwind, parameswaran2019novel, bhat2022improved} due to its ability to model fluid interface accurately, with the right jump conditions. Since the incompressible Navier-Stokes equation is a system of elliptic-parabolic partial differential equations (PDE), development of Riemann solvers necessitates a change in the nature of the system of PDEs. The most widely used approach is the artificial compressibility (AC) method \cite{chorin1967numerical} which leads to a hyperbolic-parabolic system of PDEs. When extended to unsteady flows, AC method utilizes dual time-stepping \cite{gaitonde1998dual, nithiarasu2003efficient} which involves excessive iteration in pseudo-time, to ensure a divergence-free velocity field, rendering the method computationally inefficient. Recently, weakly compressible (WC) methods have emerged as a computationally efficient alternative for solving incompressible flow problems \cite{delorme2017simple, shi2020simulations} which also admits development of Riemann solvers \cite{pan2022high}. However, when extended to two-phase flows,  WC models lead to a non-conservative system which poses a particular challenge in a finite volume framework \cite{volpert1967spaces, dal1995definition}.

In literature, non-conservative systems have been extensively studied in the context of shallow-water equations \cite{chandrashekar2020path} and compressible multiphase flows \cite{saurel1999multiphase, de2018hllc}, which can pave the way for developing hyperbolic solvers for WC two-phase systems. In particular, the finite volume treatment of compressible multi-fluid models such as the Baer-Nunziato model \cite{baer1986two} and its derivatives \cite{saurel2009simple, kapila2001two}, due to the similarities in underlying physics, motivates extension of the same to WC two-phase models. To this end, in our recent paper \cite{melvin2025development}, we had proposed a Riemann solver for the conservative fluxes coupled with a discretization for non-conservative terms which is consistent with a steady-state constraint \cite{abgrall1996prevent}, for a WC two-phase model. While this approach gives a simple oscillation-free discretization for the non-conservative terms, deriving closed form expressions for the same can be challenging. In particular, the steady-state constraint considered to derive the discretization for non-conservative terms \cite{melvin2025development}, needs to be modified when additional effects, such as surface tension terms \cite{garrick2017finite, oomar2021all}, are included in the non-conservative system.

Path-conservative methods \cite{pares2006numerical, castro2013hllc, castro2017well} offer an alternate approach to treat non-conservative hyperbolic systems. The theory developed in \cite{dal1995definition} defines weak solutions of non-conservative system which forms the basis of path-conservative schemes. According the theory in \cite{dal1995definition}, the non-conservative product is realized as a Borel measure, assuming the solution variable to be a function of bounded variation \cite{le1989shock}. Several path-conservative schemes have been developed for compressible multiphase models \cite{toumi1992weak, dumbser2010force, franquet2012runge, de2018hllc}. Path-conservative schemes \cite{dumbser2011universal, dumbser2011simple, dumbser2016new} have also been proposed for compressible multiphase model that includes surface tension effects in the non-conservative system \cite{nguyen2015path, nguyen2017numerical}. In this framework, the capillary effects are naturally taken into account by the hyperbolic solver, eliminating the need to discretize surface tension terms separately. However, determining the right family of paths remains one of the main limitations of path-conservative methods. Further, examples of convergence difficulties in path-conservative methods were given in \cite{castro2008many, abgrall2010comment}, some of which were addressed in \cite{castro2017well, chalons2017new} which the interested readers are referred to.

In the present work, a HLLC-type path-conservative scheme is developed for a WC two-phase model. The ability of HLLC solver to model fluid interface as a contact wave \cite{bhat2022improved} makes it an ideal choice to simulate two-phase flows. Unlike in our previous work \cite{melvin2025development}, the path-conservative method inherently includes the non-conservative terms, eliminating the need to specify a separate discretization for the non-conservative terms. Additionally, surface tensions terms are naturally included in the HLLC solver, as in \cite{nguyen2015path}, which leads to accurate jump conditions in pressure, velocity and phase function across the waves. By including the surface tension terms in the path-conservative solver, pressure jump obtained across the fluid interface is consistent with the Young-Laplace law. Therefore, additional interfacial pressure jump considerations, as in the ghost fluid method \cite{fedkiw1999non}, are not required to balance the forces and eliminate non-physical oscillations near the interface. Moreover, the need to include surface tension terms in the hyperbolic solver of WC two-phase model is demonstrated through numerical examples, later in the present work. 

The WC two-phase model, adopted in the present work, considers a pressure evolution equation obtained from the compressible energy equation and the evolution of fluid interface is modelled using the conservative phase-field (CPF) method \cite{chiu2011conservative}. In our previous work \cite{melvin2025development}, the interface dynamics was captured using the conservative level set (CLS) method \cite{olsson2005conservative}. However, in the context of WC models, the fidelity of the reinitialization process between subsequent physical time steps is contentious, as it modifies the conservative level set field while the pressure and velocity fields are unaltered \cite{so2011anti, kajzer2020weakly}. This issue is absent in the CPF method as the interface advection equation includes the interface regularization term. The CPF method can be interpreted as a single-step equivalent of CLS method \cite{chiu2011conservative, kajzer2020weakly} and therefore retains the excellent mass conservation property of the latter. Furthermore, the diffused nature of the interface in CPF method, as in CLS method, aids in accurate estimation of interface normals and curvatures which are crucial in surface tension dominated flows. The CPF method facilitates a simple gradient estimation technique to compute interface curvature, for an arbitrary mesh, unlike the sharp interface approaches, such as volume of fluid (VOF) method \cite{hirt1981volume}, where curvature computation requires special attention \cite{popinet2018numerical}. Additionally, material properties such as density and viscosity are related to the phase-field variable though a linear expression which simplifies the formulation of the path-conservative scheme \cite{melvin2025development}. All the equations in the WC two-phase model are solved simultaneously to prevent any lag between the solution variables during time-marching \cite{yang2014upwind, yang2020hllc}. The solver is validated against benchmark results available in literature with emphasis on surface tension dominated problems. In this work, the finite volume algorithm is developed for an arbitrary mesh and is therefore tested on structured as well as unstructured grids.

The mathematical formulation of the WC two-phase model adopted in the study is introduced in section \ref{sec:math}, followed by the numerical methods for discretizing the model in section \ref{sec:num}. Section \ref{sec:res} discusses the results and inferences from the numerical experiments. Finally, the conclusions drawn from the present work are outlined in section \ref{sec:conc}.

\section{Mathematical model}\label{sec:math}

In the weakly compressible framework, the continuity equation is replaced with a pressure evolution equation. In the present work, the following pressure evolution equation is considered
\begin{equation}\label{eq:math_pe}
	\dfrac{\partial p}{ \partial t} + \mathbf{v}\cdot \nabla p + \rho a^2 \nabla \cdot \mathbf{v} = 0
\end{equation}
where $p$ is the pressure, $\rho$ is the density and ${\mathbf{v}} \equiv (u,v)$ is the velocity vector. The speed of sound $a$ is usually considered a tunable parameter in weakly compressible (WC) models. The pressure evolution equation \eqref{eq:math_pe} is derived from the compressible energy equation neglecting the heat flux and viscous dissipation terms \cite{toutant2017general}. In the low Mach regime, it was demonstrated that the pressure advection term $\mathbf{v}\cdot \nabla p$ is superfluous and can be neglected from the WC framework \cite{toutant2017general, toutant2018numerical}. However, the present work retains the pressure advection term to preserve the hyperbolicity of the system when surface tension terms are considered in the hyperbolic solver (demonstrated later). Furthermore, in order to damp the acoustic waves, pressure diffusion term is often added to the pressure evolution equation \cite{clausen2013entropically, toutant2018numerical}, which was also considered in our previous work \cite{melvin2025development}. However, efficacy of the pressure diffusion term in damping the acoustic waves is questionable \cite{raghunathan2024non} and is therefore not considered in the present study, as in \cite{matsushita2021gas, kajzer2022weakly}.

The pressure evolution equation is combined with the momentum conservation equations

\begin{subequations}\label{eq:math_mom}
	\begin{align}
		\dfrac{\partial (\rho u)}{\partial t} + \dfrac{\partial (\rho u^2 + p)}{\partial x} + \dfrac{\partial (\rho uv)}{\partial y} &= \dfrac{\partial}{\partial x}\left( 2 \mu \dfrac{\partial u}{\partial x}\right) +  \dfrac{\partial}{\partial y}\left\{ \mu \left(\dfrac{\partial u}{\partial y} + \dfrac{\partial v}{\partial x} \right) \right\} + \rho g_x + \sigma \kappa \dfrac{\partial \psi}{\partial x} \\
		\dfrac{\partial (\rho v)}{\partial t} + \dfrac{\partial (\rho uv)}{\partial x} + \dfrac{\partial (\rho v^2 + p)}{\partial y} &= \dfrac{\partial}{\partial x}\left\{ \mu \left(\dfrac{\partial u}{\partial y} + \dfrac{\partial v}{\partial x} \right) \right\} + \dfrac{\partial}{\partial y}\left( 2 \mu \dfrac{\partial v}{\partial y}\right) + \rho g_y + \sigma \kappa \dfrac{\partial \psi}{\partial y}
	\end{align}
\end{subequations}
where $\mu$ is the dynamic viscosity. The body forces are included using the source term $\rho \mathbf{g}$ where ${\mathbf{g}} \equiv (g_x,g_y)$ is the acceleration due to gravity. Capillary effects are included using the continuum surface force (CSF) \cite{brackbill1992continuum} term $\sigma \kappa \nabla \psi$ where $\sigma$ is the surface tension coefficient and $\kappa$ is the interface curvature computed from the scalar phase field variable $\psi$ as
\begin{equation}\label{eq:math_curvature}
	\kappa = -\nabla \cdot \left( \dfrac{\nabla \psi}{|\nabla \psi|} \right)
\end{equation}

To capture the interface dynamics, a conservative phase field model \cite{chiu2011conservative} is employed 
\begin{equation}\label{eq:math_pfm}
	\dfrac{\partial \psi}{\partial t} + \dfrac{\partial (u\psi)}{\partial x} + \dfrac{\partial (v\psi)}{\partial y} = \nabla \cdot \mathbf{f}_R
\end{equation}
where $\psi$ is the phase field variable and $\mathbf{f}_R$ is the regularization term defined in equation \eqref{eq:math_regular_scls}. A hyperbolic tangent function is considered as the phase field variable $\psi$
\begin{equation}\label{eq:math_pf_func}
	\psi(x,y,t) = \dfrac{1}{2}\left\{ \tanh\left(\frac{\phi(x,y,t)}{2\varepsilon}\right) + 1\right\}
\end{equation}
where $\varepsilon$ is a mesh dependent parameter that dictates the width of the smooth transition region between the two fluids and $\phi$ is the signed distance function. For the phase field variable given by equation \eqref{eq:math_pf_func}, the fluid interface is represented by the $\psi = 0.5$ contour. The material properties such as density and viscosity are related to the phase field variable through the following expressions
\begin{equation}\label{eq:math_material}
	\rho = \rho_1 \psi + \rho_2 (1 - \psi), \quad\text{and}\quad \mu = \mu_1 \psi + \mu_2 (1 - \psi),
\end{equation}
where $(\cdot)_1$ and $(\cdot)_2$ are the material properties of fluid 1 and 2 respectively.

The governing equations \eqref{eq:math_pe}, \eqref{eq:math_mom} and \eqref{eq:math_pfm} can be written in a compact form as
\begin{equation}\label{eq:math_compact}
	\dfrac{\partial \mathbf{U}}{\partial t} + \dfrac{\partial \mathbf{F}_c}{\partial x} + \dfrac{\partial \mathbf{G}_c}{\partial y} + \mathbf{B}_x \dfrac{\partial \mathbf{U}}{\partial x} + \mathbf{B}_y \dfrac{\partial \mathbf{U}}{\partial y} = \dfrac{\partial \mathbf{F}_v}{\partial x} + \dfrac{\partial \mathbf{G}_v}{\partial y} + \mathbf{F}_g + \nabla \cdot \mathbf{F}_r
\end{equation}
where

\begin{equation*}
	\begin{split}
		\mathbf{U} = 
		\begin{bmatrix}
			p/\beta \\
			\rho u \\
			\rho v \\
			\psi
		\end{bmatrix},\quad
		\mathbf{F}_c = 
		\begin{bmatrix}
			\rho u \\
			\rho u^2 + p \\
			\rho uv \\
			u \psi
		\end{bmatrix},\quad
		\mathbf{G}_c &=
		\begin{bmatrix}
			\rho v \\
			\rho uv \\
			\rho v^2 + p \\
			v \psi
		\end{bmatrix}, \quad
		\mathbf{B}_x =
		\begin{bmatrix}
			u & 0 & 0 & -(\rho_1 - \rho_2)u \\
			0 & 0 & 0 & -\sigma \kappa \\
			0 & 0 & 0 & 0 \\
			0 & 0 & 0 & 0 
		\end{bmatrix}, \\
		\mathbf{B}_y =
		\begin{bmatrix}
			v & 0 & 0 & -(\rho_1 - \rho_2)v \\
			0 & 0 & 0 & 0 \\
			0 & 0 & 0 & -\sigma \kappa \\
			0 & 0 & 0 & 0 
		\end{bmatrix},\quad
		\mathbf{F}_v &= 
		\begin{bmatrix}
			0 \\
			2\mu \frac{\partial u}{\partial x} \\
			\mu \left( \frac{\partial u}{\partial y} + \frac{\partial v}{\partial x}\right)\\
			0
		\end{bmatrix},\quad 
		\mathbf{G}_v = 
		\begin{bmatrix}
			0 \\
			\mu \left( \frac{\partial u}{\partial y} + \frac{\partial v}{\partial x}\right)\\
			2\mu \frac{\partial v}{\partial y}\\
			0
		\end{bmatrix}, \quad \\
		\mathbf{F}_g = 
		\begin{bmatrix}
			0\\
			\rho g_x \\
			\rho g_y \\
			0
		\end{bmatrix},\quad & \text{and}\quad
		\nabla \cdot \mathbf{F}_r = 
		\begin{bmatrix}
			0 \\
			0 \\
			0 \\
			\nabla \cdot \mathbf{f}_R
		\end{bmatrix}
	\end{split}
\end{equation*}
Here $\mathbf{U}$ is the vector of conserved variables, $\mathbf{F}_c$ and $\mathbf{G}_c$ are the convective fluxes in the $x$ and $y$ directions respectively. $\mathbf{B}_x$ and $\mathbf{B}_y$ are the coefficient matrices of the non-conservative products. The viscous fluxes in the $x$ and $y$ directions are denoted by $\mathbf{F}_v$ and $\mathbf{G}_v$ respectively. Lastly, $\mathbf{F}_g$ and $\mathbf{F}_r$ are the gravitational force and phase field regularization terms, respectively. For simplicity, the speed of sound $a$ in the pressure evolution equation \eqref{eq:math_pe} has been replaced with the artificial compressibility (AC) parameter $\beta = a^2$ in equation \eqref{eq:math_compact}, which is taken as a tunable parameter \cite{melvin2025development}.

Despite being dependent on the phase field variable $\psi$ \eqref{eq:math_curvature}, interface curvature $\kappa$ is considered an independent variable in equation \eqref{eq:math_compact}. The interface curvature is assumed to be a geometrical parameter that is locally constant at a given time, which simplifies the analysis of the hyperbolic system \cite{perigaud2005compressible, nguyen2015path}.

\subsubsection*{Phase field interface-regularization}

The interface-regularization term is responsible for retaining the initial property of the phase field variable throughout the simulation. In the seminal work on conservative phase field method \cite{chiu2011conservative}, the interface-regularization term was defined as
\begin{equation}\label{eq:math_regular_chiu}
	\mathbf{f}_R = \gamma \left\{ \varepsilon \nabla \psi - \psi (1 - \psi) \mathbf{n}_\psi \right\}
\end{equation}
where $\mathbf{n}_\psi = \nabla \psi / |\nabla \psi|$ is the interface normal vector and $\gamma$ is a velocity-scale parameter. The mesh dependent parameter $\varepsilon$ is taken as
\begin{equation}\label{eq:math_regular_eps}
	\varepsilon = \dfrac{h^{1-d}}{2}
\end{equation} 
where $h$ is the average mesh size and $d \in [0,0.1]$ is a tunable parameter. An interpretation of the conservative phase field method as a single-step equivalent of the conservative level set (CLS) method \cite{olsson2005conservative, olsson2007conservative} was also given in \cite{chiu2011conservative}. Therefore, the regularization term suffers the same drawbacks as the reinitialization procedure in CLS method: a) undesired interface movement, and b) unphysical patch formation away from the interface. In literature, several fixes have been proposed for the reinitialization \cite{desjardins2008accurate, shukla2010interface, mccaslin2014localized, waclawczyk2015consistent, chiodi2017reformulation, shervani2018stabilized, parameswaran2023stable} as well as the phase field regularization \cite{chiu2019coupled, jain2022accurate, hwang2024robust} procedures. These fixes come at the expense of computational efficiency and/or conservation. In the present work, the reinitialization technique from the stabilized conservative level set (SCLS) method \cite{shervani2018stabilized} is adopted for interface-regularization. In the SCLS method, the interface normal vector is estimated as \cite{shervani2018stabilized}
\begin{equation}\label{eq:math_interNormal_scls}
	\mathbf{n}_\psi = \dfrac{\nabla \psi}{\sqrt{|\nabla \psi|^2 + \varepsilon \exp\left(-\delta \varepsilon^2|\nabla \psi|^2 \right)}}
\end{equation}
where the tunable parameter $\delta = 10$, as recommended in \cite{shervani2018stabilized}. The magnitude of the interface normal estimated by equation \eqref{eq:math_interNormal_scls}, diminishes away from the fluid interface thus alleviating the unphysical patch formation away from the interface due to ill-conditioned interface normals. The regularization term can be written as 
\begin{equation}\label{eq:math_regular_scls}
	\mathbf{f}_R = \gamma \left\{ \varepsilon (\nabla \psi \cdot \mathbf{n}_\psi) \mathbf{n}_\psi + (1 - |\mathbf{n}_\psi|^2) \varepsilon \nabla \psi - \psi (1-\psi) \mathbf{n}_\psi \right\}
\end{equation}
where the first two terms are the diffusive terms and the last term is the compressive term. The velocity-scale parameter is taken as $ \gamma =  U_\text{max}(t)$, with $U_\text{max}(t)$ being the instantaneous maximum velocity magnitude in the domain \cite{chiu2011conservative}.

\section{Numerical method}\label{sec:num}

As previously outlined, the presence of non-conservative products in the governing equation poses a particular challenge in its treatment using finite volume methods. The basis of the challenge is that the traditional definition of weak solution for conservative systems \cite{leveque2002finite} cannot be extended to non-conservative systems. Fortunately, the theory of Dal Maso, Le Floch and Murat \cite{dal1995definition} gives a rigorous definition of weak solutions for non-conservative hyperbolic systems. Therefore, it is important to establish the hyperbolicity of the present model to develop a path-conservative scheme for it.

\subsection{Hyperbolicity in time}

Ignoring the terms on the right hand side of equation \eqref{eq:math_compact}, the inviscid form of the governing equation can be written as
\begin{equation}\label{eq:num_inv_2d}
	\dfrac{\partial \mathbf{U}}{\partial t} + \dfrac{\partial \mathbf{F}_c}{\partial x} + \dfrac{\partial \mathbf{G}_c}{\partial y} + \mathbf{B}_x \dfrac{\partial \mathbf{U}}{\partial x} + \mathbf{B}_y \dfrac{\partial \mathbf{U}}{\partial y} = 0
\end{equation}
The conservative and non-conservative terms satisfy the following rotational invariance property \cite{bauerle2025rotational}
\begin{equation}\label{eq:num_rotationalInv}
	\left\{ \mathbf{F}_c(\mathbf{U})n_x + \mathbf{G}_c(\mathbf{U})n_y \right\} + \left\{ \mathbf{B}_x(\mathbf{U}, \kappa) \mathbf{U} n_x + \mathbf{B}_y(\mathbf{U}, \kappa) \mathbf{U} n_y \right\} = \mathbf{T}^{-1} \left\{\mathbf{F}_c(\hat{\mathbf{U}}) + \,\mathbf{B}_x(\hat{\mathbf{U}}, \kappa) \hat{\mathbf{U}} \right\}
\end{equation}
where $\hat{\mathbf{U}} = \mathbf{T}\mathbf{U} = \left[ p/\beta, \, \rho\hat{u}, \, \rho\hat{v}, \, \psi \right]^T$ is the transformed vector of conserved variables for a given rotation matrix $\mathbf{T}$ and its inverse $\mathbf{T}^{-1}$ defined as
\begin{equation}\label{eq:num_rotationmatrix}
	\mathbf{T} = 
	\begin{bmatrix}
		1 & 0 & 0 & 0 \\
		0 & n_x & n_y & 0 \\
		0 & -n_y & n_x & 0 \\
		0 & 0 & 0 & 1
	\end{bmatrix}\quad\text{and}\quad
	\mathbf{T}^{-1} = 
	\begin{bmatrix}
		1 & 0 & 0 & 0 \\
		0 & n_x & -n_y & 0 \\
		0 & n_y & n_x & 0 \\
		0 & 0 & 0 & 1
	\end{bmatrix}
\end{equation}
In a finite volume framework, these rotations correspond to the outward facing normals $\mathbf{n} \equiv (n_x, n_y)$ of the boundaries of control volume (see \cite{melvin2025development} for details). Using the property \eqref{eq:num_rotationalInv}, the two-dimensional problem \eqref{eq:num_inv_2d} can be reduced to an augmented one-dimensional one as
\begin{equation}\label{eq:num_inv_1d}
	\dfrac{\partial \hat{\mathbf{U}}}{\partial t} + \dfrac{\partial \mathbf{F}_c(\hat{\mathbf{U}})}{\partial \hat{x}} + \mathbf{B}_{\hat{x}}(\hat{\mathbf{U}}, \kappa) \dfrac{\partial \hat{\mathbf{U}}}{\partial \hat{x}} = 0
\end{equation}
The above equation \eqref{eq:num_inv_1d} can be written in the following quasi-linear form
\begin{equation}\label{eq:num_inv_quasiLinear}
	\dfrac{\partial \hat{\mathbf{U}}}{\partial t} + \mathbf{A}(\hat{\mathbf{U}}, \kappa) \dfrac{\partial \hat{\mathbf{U}}}{\partial \hat{x}} = 0
\end{equation}
where the coefficient matrix $\mathbf{A}(\hat{\mathbf{U}})$ is defined as

\begin{equation*}
	\mathbf{A}(\hat{\mathbf{U}}, \kappa) = \dfrac{\partial \mathbf{F}_c(\hat{\mathbf{U}})}{\partial \hat{\mathbf{U}}} + \mathbf{B}_{\hat{x}}(\hat{\mathbf{U}}, \kappa) = 
	\begin{bmatrix}
		\hat{u} & 1 & 0 & -(\rho_1 - \rho_2)\hat{u} \\
		\beta & 2\hat{u} & 0 & -(\rho_1 - \rho_2)\hat{u}^2 - \sigma \kappa \\
		0 & \hat{v} & \hat{u} & -(\rho_1 - \rho_2)\hat{u}\hat{v} \\
		0 & \dfrac{\psi}{\rho} & 0 & \dfrac{\rho_2 \hat{u}}{\rho}
	\end{bmatrix}
\end{equation*}
The coefficient matrix $\mathbf{A}$ has four real eigenvalues
\begin{equation}\label{eq:num_eigenvalues}
	\lambda_1 = \hat{u}_\rho - a_s,\quad \lambda_2 = \lambda_3 = \hat{u}\text{ and }\lambda_4 = \hat{u}_\rho + a_s
\end{equation}
where

\[ \hat{u}_\rho = \left(1 + \dfrac{\rho_2}{2\rho}\right) \hat{u} \quad\text{and}\quad a_s = \sqrt{(\hat{u}_\rho - \hat{u})^2 + \beta - \dfrac{\sigma \kappa \psi}{\rho}} \]
and the corresponding right eigenvector matrix, written as

\begin{equation*}\label{eq:num_eigenvectors}
	\mathbf{R} = 
	\begin{bmatrix}
		\mathbf{R}_1 & \mathbf{R}_2 & \mathbf{R}_3 & \mathbf{R}_4
	\end{bmatrix} = 
	\begin{bmatrix}
		1 & \dfrac{\sigma \kappa}{\beta} & 0 & 1 \\
		\left( \lambda_1 - \dfrac{\rho_2 \hat{u}}{\rho} \right) & (\rho_1 - \rho_2)\hat{u} & 0 & \left( \lambda_4 - \dfrac{\rho_2 \hat{u}}{\rho} \right) \\
		\hat{v} & 0 & 1 & \hat{v} \\
		\dfrac{\psi}{\rho} & 1 & 0 & \dfrac{\psi}{\rho}
	\end{bmatrix}
\end{equation*}
with four linearly independent columns, thus verifying that the inviscid subsystem is hyperbolic in time.

\begin{remark}
	The inclusion of surface tension in the non-conservative system modifies the eigenvalues in the present model, unlike in the compressible multiphase models \cite{nguyen2017numerical}. It can be attributed to conservative advection equation \eqref{eq:math_pfm} adopted in the present work. However, employing a non-conservative advection equation
	
	\[ \dfrac{\partial \psi}{\partial t} + \mathbf{v} \cdot \nabla \psi = \nabla \cdot \mathbf{f}_R \]
	instead of equation \eqref{eq:math_pfm} will lead to a non-conservative system with its eigenvalues independent of surface tension terms, as in \cite{nguyen2017numerical}. The analysis is beyond the scope of this work.
\end{remark}

\subsection{Generalized Riemann invariant analysis}

A generalized Riemann invariant analysis \cite{toro2013riemann} of the intermediate waves namely, the contact and shear waves yield the accurate jump conditions which are essential in obtaining the intermediate states of the Riemann solver. The generalized Riemann invariants across the contact wave $\mathbf{R}_2$ can be written as
\begin{equation}\label{eq:num_gri_contact}
	\dfrac{d(p/\beta)}{\sigma \kappa / \beta} = \dfrac{d(\rho \hat{u})}{(\rho_1 - \rho_2)\hat{u}} = \dfrac{d(\rho \hat{v})}{0} = \dfrac{d \psi}{1}
\end{equation}
and across the shear wave $\mathbf{R}_3$ can be written as
\begin{equation}\label{eq:num_gri_shear}
	\dfrac{d(p/\beta)}{0} = \dfrac{d(\rho \hat{u})}{0} = \dfrac{d(\rho \hat{v})}{1} = \dfrac{d \psi}{0}
\end{equation}

Across the contact wave $\mathbf{R}_2$, through mathematical manipulation it can be shown that the normal velocity $\hat{u}$ remains constant while tangential velocity $\hat{v}$ and phase field variable $\psi$ vary. Pressure jump across the contact wave is obtained as
\begin{equation}\label{eq:num_laplacelaw}
	p^+ - p^- = \sigma \kappa (\psi^+ - \psi^-)
\end{equation}
which is consistent with the Young-Laplace law. Here $(\cdot)^-$ and $(\cdot)^+$ are the limits of the variables to left and right of the contact wave. In contrast to our previous work \cite{melvin2025development}, where pressure was constant across the contact wave, inclusion of surface tension terms in the non-conservative system yields a jump in pressure. Lastly, across the shear wave $\mathbf{R}_3$, only the tangential velocity varies while the other variables remain constant, as in \cite{melvin2025development}.

\subsection{Path-conservative method}

Being hyperbolic in time, the present non-conservative system admits development of path-conservative finite volume schemes. Here, the theory necessary for the development of path-conservative schemes for non-conservative systems is briefly reviewed.

Casting the governing equation in integral form forms the basis for developing finite volume schemes. However, when discontinuities are present in the solution, integral formulation requires a definition of weak solution \cite{leveque2002finite}. In standard finite volume techniques, a weak solution is defined by multiplying equation \eqref{eq:num_inv_quasiLinear} with a smooth test function (with compact support) and transferring the spatial derivatives to it. Thus, the discontinuities present in the solution variable no longer pose a challenge when integrated over space. However, it is essential that the coefficient matrix $\mathbf{A}$ be the Jacobian of the flux function for this definition of weak solution and therefore, fails for non-conservative systems such as equation \eqref{eq:num_inv_1d}. Weak solution for the non-conservative system is given by the theory developed in \cite{dal1995definition} where the non-conservative product is realized as a Borel measure. Defining this locally bounded measure requires a choice of family of Lipschitz continuous paths $\Phi: [0,1] \times \Omega \times \Omega \mapsto \Omega$ connecting the left $\mathbf{U}_L$ and right $\mathbf{U}_R$ states in a set of all admissible states $\Omega$. The family of paths must satisfy the following conditions

\[ \Phi(0; \mathbf{U}_L,\mathbf{U}_R) = \mathbf{U}_L, \quad \Phi(1; \mathbf{U}_L,\mathbf{U}_R) = \mathbf{U}_R\]
and

\[ \Phi(s; \mathbf{U},\mathbf{U}) = \mathbf{U}\]
In practice, the family of paths give a sense to the integral of the non-conservative product, over space \cite{castro2017well}. Denoting this spatial integral over an interval $[x_1,x_2]$ as

\[ \fint_{x_1}^{x_2} \mathbf{A}(\mathbf{U}) \dfrac{\partial \mathbf{U}}{\partial x}\, dx\text{,} \]
its Lebesgue decomposition \cite{hewitt2012real} can be written as
\begin{equation}\label{eq:num_lebesgue}
	\mu_a^\Phi = \int_{x_1}^{x_2} \mathbf{A}(\mathbf{U}) \dfrac{\partial \mathbf{U}}{\partial x}\, dx \quad \text{and} \quad
	\mu_s^\Phi = \sum_{l} \left( \int_0^1 \mathbf{A} (\Phi(s;\mathbf{U}^-_l,\mathbf{U}^+_l)) \dfrac{\partial \Phi}{\partial s} \, ds \right)\delta_{x=x_l}
\end{equation}
where $\mu_a^\Phi$ and $\mu_s^\Phi$ are the absolutely continuous and singular continuous parts respectively. Here $\delta_{x=x_l}$ is the Dirac delta function at the $l^\text{th}$ discontinuity located at $x_l \in [x_1,x_2]$, and $\mathbf{U}^-_l$ and $\mathbf{U}^+_l$ are the limits of the solution variable $\mathbf{U}$ to the left and right of the $l^\text{th}$ discontinuity.

In a first-order finite volume scheme, where the solution is piecewise constant within a control volume, the measure only consists of the singular part with the punctual masses $\mu_s^\Phi$ placed at the cell interfaces. The Riemann problem at cell interfaces has the same configuration, as in the conservative systems, consisting of discontinuities: contact, shock, and rarefaction waves. Across a discontinuity with speed $\xi$, a weak solution should satisfy the generalized Rankine-Hugoniot jump condition \cite{pares2006numerical}
\begin{equation}\label{eq:num_generalizedRH}
	\int_0^1 \left\{ \mathbf{A} (\Phi(s;\mathbf{U}^-,\mathbf{U}^+)) - \xi \mathbf{I} \right\} \dfrac{\partial \Phi}{\partial s} \, ds = 0
\end{equation}
where $\mathbf{I}$ is the identity matrix and $\mathbf{U}^-$, $\mathbf{U}^+$ are the solution states to the left and right of the discontinuity. For a conservative system, it can be shown that the generalized Rankine-Hugoniot condition \eqref{eq:num_generalizedRH} reduces to the usual Rankine-Hugoniot condition \cite{leveque2002finite}.

Path-conservative schemes aim to split the punctual masses, often called fluctuations, at the cell interfaces to identify its contribution to the cells left and right of the interface. The present work aims to develop a HLLC-type path-conservative scheme to split these fluctuations. Interested readers are referred to \cite{pares2006numerical, castro2017well} for a detailed review of various path-conservative schemes.

\subsection{First-order HLLC-type path-conservative scheme}

Based on the theory reviewed in the last section, the present work aims to develop a HLLC-type path-conservative scheme to split the fluctuations at the cell interfaces, for the inviscid subsystem \eqref{eq:num_inv_quasiLinear}. 

\begin{figure}[ht!]
	\centering
	\includegraphics[width=4in]{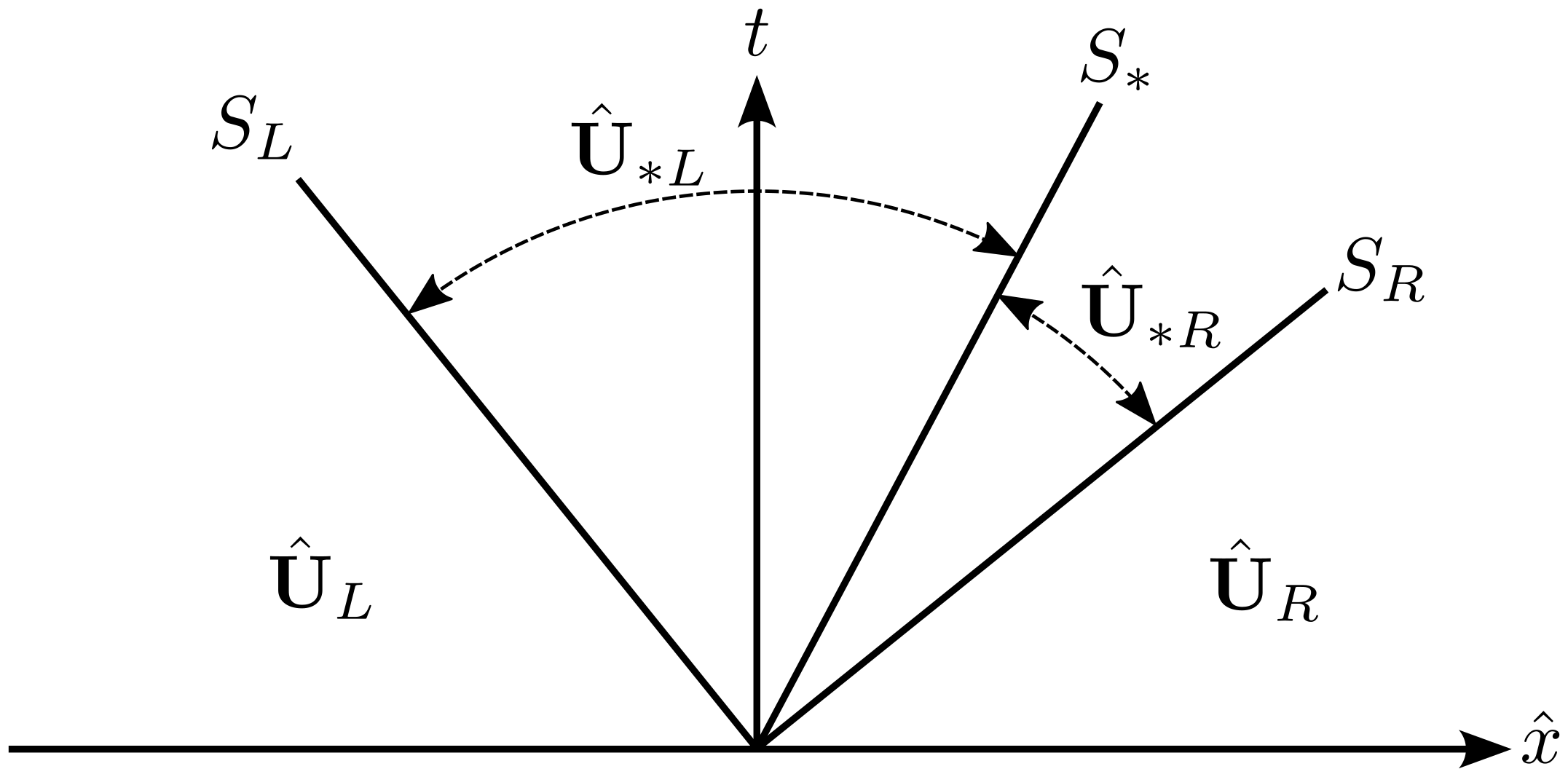}
	\caption{The wave structure of the HLLC Riemann solver.}
	\label{fig:num_hllc_waves}
\end{figure}

To develop a HLLC solver for the system, a three wave structure, as in figure \ref{fig:num_hllc_waves}, is assumed. Here, the intermediate wave models the contact and shear waves. Across each wave, the generalized Rankine-Hugoniot jump conditions \eqref{eq:num_generalizedRH} must be satisfied. However, this requires a choice of family of paths $\Phi$. A good choice for family of paths is based on the viscous profiles obtained by considering the vanishing viscosity solutions of equation \eqref{eq:num_inv_quasiLinear} \cite{le1989shock}. However, computing these viscous profiles can be a challenging task and in such cases, family of straight line segments prove a sensible choice as their jump conditions give a third-order approximation of the physically correct ones \cite{castro2013hllc}. In the present work, the canonical path
\begin{equation}\label{eq:num_canon_path}
	\Phi(s; \mathbf{U}^-,\mathbf{U}^+) = \mathbf{U}^- + s(\mathbf{U}^+ - \mathbf{U}^-) 
\end{equation}
connecting the states $\mathbf{U}^-$ and $\mathbf{U}^+$, is considered. Therefore, the generalized Rankine-Hugoniot condition \eqref{eq:num_generalizedRH} simplifies as
\begin{equation}\label{eq:num_RH_canon}
	\mathbf{F}_c(\mathbf{U}^+) - \mathbf{F}_c(\mathbf{U}^-) + \mathbf{B}_x(\mathbf{U}^\pm, \kappa)(\mathbf{U}^+ - \mathbf{U}^-) = \xi(\mathbf{U}^+ - \mathbf{U}^-)
\end{equation}
where
\begin{equation}\label{eq:num_noncons_matrixInt}
	\mathbf{B}_x(\mathbf{U}^\pm, \kappa) = \int_0^1 \mathbf{B}_x(\Phi(s;\mathbf{U}^-,\mathbf{U}^+), \kappa) \, ds
\end{equation}
In literature, several approaches to evaluate the above integral, analytically and numerically \cite{dumbser2011universal, dumbser2011simple, kamath2009roe}, have been proposed. In the present work, the integral is evaluated analytically. It is worth reiterating that the interface curvature $\kappa$ is considered locally constant and is not solved by the Riemann solver.

To obtain the intermediate states $\hat{\mathbf{U}}_{*L}$ and $\hat{\mathbf{U}}_{*R}$ in the three-wave structure (figure \ref{fig:num_hllc_waves}), the generalized Rankine-Hugoniot conditions \eqref{eq:num_RH_canon} are applied across the three waves leading to an underdetermined system of equations. However, the constraints obtained from the generalized Riemann invariant analysis completes the system, yielding the following expressions for the intermediate states:

\begin{subequations}\label{eq:num_hllcStar}
	\begin{align}
		\intertext{Pressures}
		\begin{split}
			& \left(\dfrac{p}{\beta}\right)_{*L} = \left(\dfrac{p}{\beta}\right)_{L} + \dfrac{ \rho_{*L}(S_* - \langle \hat{u} \rangle_L) - \rho_{L}(\hat{u}_L - \langle \hat{u} \rangle_L)}{S_L - \langle \hat{u} \rangle_L}, \quad \text{and} \\
			& \left(\dfrac{p}{\beta}\right)_{*R} = \left(\dfrac{p}{\beta}\right)_{R} + \dfrac{ \rho_{*R}(S_* - \langle \hat{u} \rangle_R) - \rho_{R}(\hat{u}_R - \langle \hat{u} \rangle_R)}{S_R - \langle \hat{u} \rangle_R} \\
		\end{split}
		\intertext{Normal velocity which is also considered to be the speed of contact wave}
		\hat{u}_{*L}, \hat{u}_{*R} = \hat{u}_* =& \dfrac{S_L(\rho \hat{u})_L - S_R(\rho \hat{u})_R + (\rho \hat{u}^2+p)_R - (\rho \hat{u}^2+p)_L - \sigma \kappa (\psi_R - \psi_L)}{S_L\rho_L - S_R\rho_R + (\rho_1 - \rho_2)\left\{(\hat{u}\psi)_R - (\hat{u}\psi)_L\right\}} \equiv S_* \\
		\intertext{Tangential velocities}
		(\rho \hat{v})_{*L} &= \dfrac{S_L(\rho \hat{v})_L - (\rho \hat{u}\hat{v})_L}{S_L - S_*} \quad \text{and} \quad (\rho \hat{v})_{*R} = \dfrac{S_R(\rho \hat{v})_R - (\rho \hat{u}\hat{v})_R}{S_R - S_*} \\
		\intertext{Phase-field variables}
		&\psi_{*L} = \dfrac{S_L \psi_L - (\hat{u}\psi)_L}{S_L - S_*} \quad \text{and} \quad \psi_{*R} = \dfrac{S_R \psi_R - (\hat{u}\psi)_R}{S_R - S_*}
	\end{align}
\end{subequations}
Here, $ \langle \hat{u} \rangle_L$ and $ \langle \hat{u} \rangle_R$ result from the integral \eqref{eq:num_noncons_matrixInt} in the generalized Rankine-Hugoniot jump conditions. The left and right wave speeds $S_L$ and $S_R$ are estimated as in \cite{davis1988simplified}

\begin{equation*}
	S_L = \min\{(\lambda_1)_L,(\lambda_1)_R\}\quad\text{and}\quad S_R = \max\{(\lambda_4)_L,(\lambda_4)_R\}
\end{equation*}
where $\lambda_1$ and $\lambda_4$ are as defined in equation \eqref{eq:num_eigenvalues}.

The obtained intermediate states aid in splitting the fluctuations at the cell interfaces. At a cell interface with the states $\mathbf{U}_L$ and $\mathbf{U}_R$ to its left and right respectively, the total fluctuation $\mathbf{D}(\mathbf{U}_L,\mathbf{U}_R)$ is split as

\begin{equation}
	\mathbf{D}(\mathbf{U}_L,\mathbf{U}_R, \kappa_{LR}) = \mathbf{D}^-(\mathbf{U}_L,\mathbf{U}_R, \kappa_{LR}) + \mathbf{D}^+(\mathbf{U}_L,\mathbf{U}_R, \kappa_{LR})
\end{equation}
Here $\mathbf{D}^-$ and $\mathbf{D}^+$ are the contributions to the left and right cells, respectively and $\kappa_{LR}$ is the interface curvature computed at the cell interface. For the HLLC wave structure considered in figure \ref{fig:num_hllc_waves}, these contributions are defined as \cite{castro2017well}:

\begin{subequations}\label{eq:num_hllc_fluctuations}
	\begin{align}
		\mathbf{D}^-(\hat{\mathbf{U}}_L, \hat{\mathbf{U}}_R, \kappa_{LR}) & =
		\begin{cases}
			0, & S_L \geq 0 \\
			S_L(\hat{\mathbf{U}}_{*L} - \hat{\mathbf{U}}_L), & S_L < 0 \leq S_* \\
			S_L(\hat{\mathbf{U}}_{*L} - \hat{\mathbf{U}}_L) + S_*(\hat{\mathbf{U}}_{*R} - \hat{\mathbf{U}}_{*L}), & S_* < 0 < S_R \\
			S_L(\hat{\mathbf{U}}_{*L} - \hat{\mathbf{U}}_L) + S_*(\hat{\mathbf{U}}_{*R} - \hat{\mathbf{U}}_{*L}) + S_R(\hat{\mathbf{U}}_{R} - \hat{\mathbf{U}}_{*R}), & S_R \leq 0
		\end{cases} \\
		\intertext{and}
		\mathbf{D}^+(\hat{\mathbf{U}}_L, \hat{\mathbf{U}}_R, \kappa_{LR}) & =
		\begin{cases}
			S_L(\hat{\mathbf{U}}_{*L} - \hat{\mathbf{U}}_L) + S_*(\hat{\mathbf{U}}_{*R} - \hat{\mathbf{U}}_{*L}) + S_R(\hat{\mathbf{U}}_{R} - \hat{\mathbf{U}}_{*R}), & S_L \geq 0 \\
			S_*(\hat{\mathbf{U}}_{*R} - \hat{\mathbf{U}}_{*L}) + S_R(\hat{\mathbf{U}}_{R} - \hat{\mathbf{U}}_{*R}), & S_L < 0 \leq S_* \\
			S_R(\hat{\mathbf{U}}_{R} - \hat{\mathbf{U}}_{*R}), & S_* < 0 < S_R \\
			0, & S_R \leq 0
		\end{cases}
	\end{align}
\end{subequations}

Using the HLLC scheme, the one-dimensional semi-discrete form of the quasi-linear equation \eqref{eq:num_inv_quasiLinear} can be written as
\begin{equation}\label{eq:num_hllc_semidiscrete}
	\dfrac{d \mathbf{U}_i}{d t} + \dfrac{1}{\Delta x} \left\{\mathbf{D}^+(\mathbf{U}_{i-1},\mathbf{U}_{i}, \kappa_{i-1/2}) + \mathbf{D}^-(\mathbf{U}_{i},\mathbf{U}_{i+1}, \kappa_{i+1/2}) \right\} = 0
\end{equation}
where $\Delta x$ is the space step. The two-dimensional form of the semi-discrete equation \eqref{eq:num_hllc_semidiscrete}, along with its high-order extension, is presented later.

\subsection{High-order extension based on solution reconstruction}

Similar to the conservative hyperbolic problems, higher-order spatial accuracy for non-conservative systems can be obtained using a solution reconstruction technique. Unlike the first-order scheme, high-order extension requires the computation of the absolute continuous part in the Lebesgue decomposition \eqref{eq:num_lebesgue}. In two-dimension, the continuous part is represented by an area integral over a control volume $\Omega$
\begin{equation}\label{eq:num_highOrder_areaInt}
	\mu_a = \iint\displaylimits_\Omega \left( \mathbf{A}_x(\mathbf{U}, \kappa) \dfrac{\partial \mathbf{U}}{\partial x} + \mathbf{A}_y(\mathbf{U}, \kappa) \dfrac{\partial \mathbf{U}}{\partial y} \right)\,d\Omega
\end{equation}
where 

\[ \mathbf{A}_x(\mathbf{U}, \kappa) = \dfrac{\partial \mathbf{F}_c (\mathbf{U})}{\partial \mathbf{U}} + \mathbf{B}_x(\mathbf{U}, \kappa) \quad\text{and}\quad \mathbf{A}_y(\mathbf{U}, \kappa) = \dfrac{\partial \mathbf{G}_c (\mathbf{U})}{\partial \mathbf{U}} + \mathbf{B}_y(\mathbf{U}, \kappa) \]
It is common to estimate this area integral by numerical integration using trapezoidal rule or Romberg extrapolation \cite{noelle2006well}. For the present model, by considering linear reconstruction within a control volume, the area integral can be simplified as
\begin{equation}\label{eq:num_highOrder}
	\mu_a = \oint\displaylimits_\Gamma (\mathbf{F}_c n^\Gamma_x + \mathbf{G}_c n^\Gamma_y) \,d\Gamma + \left\{ \mathbf{B}_x(\overline{\mathbf{U}}, \kappa_\Omega) \left. \dfrac{\partial \mathbf{U}}{\partial x}\right|_\Omega + \mathbf{B}_y(\overline{\mathbf{U}}, \kappa_\Omega) \left. \dfrac{\partial \mathbf{U}}{\partial y}\right|_\Omega \right\} \Omega
\end{equation}
where $d\Gamma$ is the infinitesimal length over the boundary $\Gamma$ of the control volume $\Omega$ and $\mathbf{n}^\Gamma \equiv (n^\Gamma_x,n^\Gamma_y)$ its outward facing unit normal vector. $\overline{\mathbf{U}}$ is the cell averaged conserved variable vector and $\kappa_\Omega$ is the interface curvature computed at the cell centre. The cell centred gradient $\left. \nabla \mathbf{U} \right|_\Omega$ is obtained from the linear reconstruction procedure. In the linear reconstruction procedure, a weighted least square technique procedure \cite{barth1989design} with geometric weights is used to obtain the gradients at the cell centres. For stability, Moore's neighbors \cite{packard1985two} are considered while formulating the least square problem. To prevent Gibbs-type oscillations, the limited gradients  $\left. \nabla \mathbf{U} \right|_\Omega$ are used, such that no new local extrema is formed at the cell interfaces during the reconstruction procedure \cite{barth1989design}.

The high-order HLLC path-conservative scheme is obtained by combining the first-order fluctuation splitting \eqref{eq:num_hllc_fluctuations} and the continuous part \eqref{eq:num_highOrder}.

\begin{remark}
	For the present model, the simplification \eqref{eq:num_highOrder} of the area integral \eqref{eq:num_highOrder_areaInt} is possible for the specific case of linear reconstruction within the control volumes. For quadratic and higher-order polynomial reconstruction, it is necessary to compute the area integral numerically \cite{noelle2006well}.
\end{remark}

\subsection{Finite volume discretization}

\begin{figure}[ht!]
	\centering
	\includegraphics[width=0.6\textwidth]{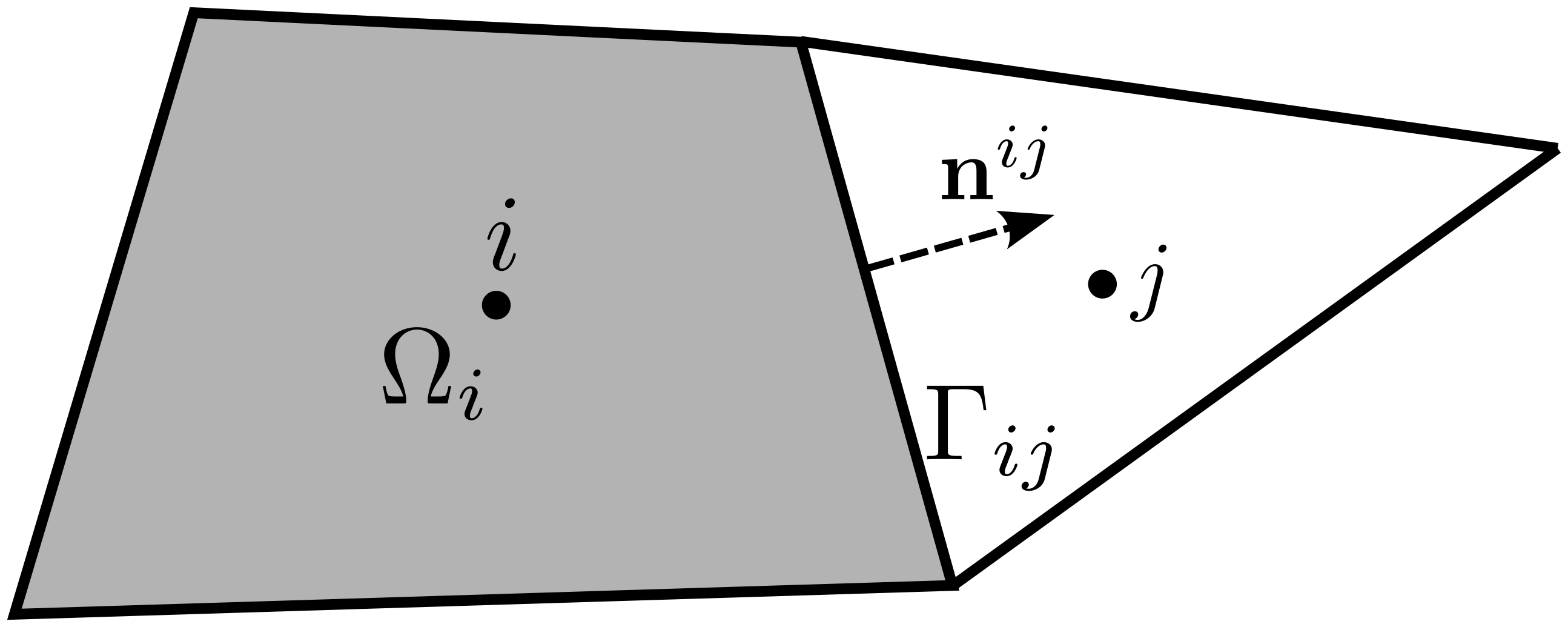}
	\caption{The $i^\text{th}$ finite volume cell. Here $\Omega_i$ denotes the area of the cell and $j$ is one of the von Neumann neighbours of the $i^\text{th}$ cell. $\mathbf{n}^{ij}$ is unit normal vector to the interface between cells $i$ and $j$, whose edge length is $\Gamma_{ij}$.}
	\label{fig:num_fvmCell}
\end{figure}

In a finite volume framework, the computational domain is discretized into a finite number of non-overlapping cells. For the $i^\text{th}$ cell (see figure \ref{fig:num_fvmCell}), spatially discretized form of equation \eqref{eq:math_compact} can be written as
\begin{equation}\label{eq:num_spaceDisc}
	\Omega_i \dfrac{\partial \overline{\mathbf{U}}_i}{\partial t} + \mathcal{R}_i(\overline{\mathbf{U}}) = 0
\end{equation}
where the residual $\mathcal{R}$ is defined as

\begin{equation*}\label{eq:num_residual}
	\begin{split}
		\mathcal{R}_i(\overline{\mathbf{U}}) &= \sum_{j \in \mathcal{V}(i)} \mathbf{T}^{-1}_{ij} \left\{\mathbf{F}_c (\hat{\mathbf{U}}_i^{ij}) + \mathbf{D}^-(\hat{\mathbf{U}}_i^{ij}, \hat{\mathbf{U}}_j^{ij}, \kappa_{ij}) \right\} \Gamma_{ij} \\ &+ \left\{ \mathbf{B}_x(\overline{\mathbf{U}}_i, \kappa_{i})\left. \dfrac{\partial \mathbf{U}}{\partial x}\right|_i + \mathbf{B}_y(\overline{\mathbf{U}}_i, \kappa_{i}) \left. \dfrac{\partial \mathbf{U}}{\partial y}\right|_i \right\} \Omega_i \\ &- \sum_{j \in \mathcal{V}(i)} (\mathbf{F}^{ij}_v n^{ij}_x + \mathbf{G}^{ij}_v n^{ij}_y) \Gamma_{ij} - \Omega_i \mathbf{F}_g(\overline{\mathbf{U}}_i) - \sum_{j \in \mathcal{V}(i)} \mathbf{F}_r^{ij} \cdot \mathbf{n}^{ij} \Gamma_{ij}
	\end{split}
\end{equation*}
Here $\overline{\mathbf{U}}_i$ is the cell averaged conserved variable vector for the $i^\text{th}$ cell whose area is denoted by $\Omega_i$. $\hat{\mathbf{U}}$ is the rotated conserved variable vector as in equation \eqref{eq:num_rotationalInv}. The rotational matrix $\mathbf{T}_{ij}$ and its inverse $\mathbf{T}_{ij}^{-1}$ \eqref{eq:num_rotationmatrix} are constructed with the unit normal vector $\mathbf{n}^{ij} \equiv (n^{ij}_x,n^{ij}_y)$ to the interface between $i^\text{th}$ and $j^\text{th}$ cells (see figure \ref{fig:num_fvmCell}) whose length is $\Gamma_{ij}$. Here $j$ is a von Neumann neighbour \cite{packard1985two} of the $i^\text{th}$ cell; the set of von Neumann neighbours is denoted by $\mathcal{V}(i)$.

The conserved variable vector $\hat{\mathbf{U}}_i^{ij}$ is the reconstructed state on the interface ${ij}$, based on the cell averaged $\overline{\mathbf{U}}_i$ and cell-centred \emph{limited} gradient $\nabla \mathbf{U}|_i$ obtained from the reconstruction procedure. Similarly,  $\hat{\mathbf{U}}_j^{ij}$ is the reconstructed state on the interface ${ij}$, based on $\overline{\mathbf{U}}_j$ and $\nabla \mathbf{U}|_j$. The fluctuation $\mathbf{D}^-(\hat{\mathbf{U}}_i^{ij}, \hat{\mathbf{U}}_j^{ij}, \kappa_{ij})$ is calculated using the HLLC fluctuation splitting scheme \eqref{eq:num_hllc_fluctuations}. Geometrical variable $\kappa_{ij}$ is the interface curvatures computed at the ${ij}^\text{th}$ face while $\kappa_{i}$ is the same computed at the $i^\text{th}$ cell centre.

The viscous fluxes $\mathbf{F}_v$ and $\mathbf{G}_v$ are computed using Green-Gauss approach along a Coirier diamond path around the cell interfaces \cite{coirier1994adaptively}. The gravitational source term is computed as the product of the cell area $\Omega_i$ and the cell averaged value of gravitational force $\mathbf{F}_g$. Discretization of the regularization term $\mathbf{F}_r$ follows the same procedure as detailed in \cite{melvin2025development} for the reinitialization equation.

\subsection{Computation of interface curvature}

Accurate computation of interface curvature \eqref{eq:math_curvature} is crucial in the present solver. As detailed previously, the first-order path-conservative scheme \eqref{eq:num_hllc_fluctuations} requires the computation of curvature at cell interfaces and the higher-order extension \eqref{eq:num_highOrder} requires the same at the cell centres. In the present work, the interface curvatures at the cell centres are computed using a weighted least square (WLS) technique \cite{barth1989design} as:

\begin{enumerate}[label=\roman*.]
	\item Through a WLS procedure, the gradients of phase-field variables $\nabla \psi$ at the cell centres and in turn the interface normals $\mathbf{n}_\psi \equiv (n_{\psi,x},n_{\psi,y}) = \nabla \psi/ |\nabla \psi|$ are computed.
	
	\item The WLS procedure is repeated using the computed interface normals to obtain the cell-centred gradients of two components of interface normals $\nabla n_{\psi,x}$ and $\nabla n_{\psi,y}$.
	
	\item Finally, the cell-centred interface curvature is obtained as
	
	\[ \kappa = - \left\{ \dfrac{\partial \left( n_{\psi,x} \right)}{\partial x} + \dfrac{\partial \left( n_{\psi,y} \right)}{\partial y} \right\} \]

\end{enumerate}

The WLS technique uses distance-based weights \cite{barth1989design} and the overdetermined system is solved through a QR factorization that utilizes Householder transformation. The interface curvatures at cell interfaces are approximated as the average of the two cell-centred curvatures either side of it \cite{garrick2017finite}.

\subsection{Temporal discretization}

A two stage strong stability preserving Runge-Kutta (SSP-RK) method \cite{gottlieb2009high} is used to discretize the time derivatives. The solution variable at $n^\text{th}$ time level $\overline{\mathbf{U}}^n$ is updated to $\overline{\mathbf{U}}^{n+1}$ for the $i^\text{th}$ cell as
\begin{equation}\label{eq:num_ssprk3}
	\begin{split}
		\overline{\mathbf{U}}^{(1)}_i &= \overline{\mathbf{U}}^n_i - \dfrac{\Delta t}{\Omega_i}\mathcal{R}_i\left( \overline{\mathbf{U}}^n \right) \\
		\overline{\mathbf{U}}^{n+1}_i &= \dfrac{1}{2}\overline{\mathbf{U}}^n_i + \dfrac{1}{2}\overline{\mathbf{U}}^{(1)}_i - \dfrac{1}{2}\dfrac{\Delta t}{\Omega_i}\mathcal{R}_i\left( \overline{\mathbf{U}}^{(1)} \right)
	\end{split}
\end{equation}
where $\overline{\mathbf{U}}^{(1)}$ is the solution variable at the intermediate stage. The global time-step $\Delta t$ is restricted by the stability requirements of the scheme. The local time-step of the $i^\text{th}$ cell $\Delta t_i$ is dictated by the maximum allowable time-steps due to convective, viscous, gravitational and surface tension terms \cite{melvin2025development}. The global time-step is taken as
\begin{equation}\label{eq:num_timeStep}
	\Delta t = \text{CFL}\,\min_{i}\{\Delta t_i\} 
\end{equation}
where the Courant number $\text{CFL} < 1$ for stability.

\section{Numerical results and discussions}\label{sec:res}

The efficacy of the developed path-conservative scheme is tested on some two-dimensional problems. The incompressible two-phase problems are chosen such that various components of the weakly-compressible solver are tested. Some of the problems are tested on structured as well as unstructured meshes to demonstrate the solver's adaptability. The results from the simulations are compared against analytical and numerical results reported in literature.

The problems are tested using first-order path-conservative scheme as well as with linear reconstruction. The \emph{effective} spatial order of convergence reported were obtained using a linear regression fit. For stability, as recommended in \cite{olsson2007conservative}, $d = 0.1$ is considered for the mesh dependent parameter $\varepsilon$ \eqref{eq:math_regular_eps}, which dictates the width of the phase-field transition region. To mimic incompressible flow, the artificial compressibility (AC) parameter $\beta$ is set appropriately high such that the Mach number is less than $0.1$ \cite{melvin2025development}. The Courant number is taken as $\text{CFL} = 0.9$ in all the problems.

\subsection{Static drop}

The case of static drop in equilibrium is an apt test case to validate the path-conservative scheme developed in the present work. The two-dimensional static drop case analysed in \cite{francois2006balanced} is examined in the present work. The problem consists of a $L \times L$ domain where $L = 8$ m with a drop of radius $R = 2$ m located at its centre. The density inside and outside the drop is considered the same $\rho_1 = \rho_2 = 1$ kg/m$^3$ and the surface tension coefficient is taken as $\sigma = 73$ N/m. The exact pressure jump across the drop is given by the Young-Laplace law

\[ \Delta p_\text{exact} = \sigma \kappa = 36.5 \text{ Pa} \]
where the interface curvature $\kappa = 1/R$.

The case of static drop is simulated using the first-order path-conservative scheme developed in this work. The numerical simulations are carried out on four uniform Cartesian meshes of sizes: $16 \times 16$, $32 \times 32$, $64 \times 64$ and $128 \times 128$. The AC parameter is taken as $\beta = 1\times 10^4$. The phase field regularization term \eqref{eq:math_regular_scls} is not considered since the goal of the test case is to validate the path-conservative scheme. Furthermore, minimal interface diffusion is expected in this static problem which eliminates the need for phase field regularization. Slip-wall boundary conditions are imposed on all four boundaries and the simulation was run till $t = 5$ s. The initial pressure and velocity is set to zero. 

The $L_2$ errors in pressure field, with respect to the exact jump $\Delta p_\text{exact}$, is estimated as

\[ L_2(p) =  \dfrac{1}{\Delta p_\text{exact}} \sqrt{ \sum_{i} \left(p_i - p^\text{exact}_i \right)^2 \,  \Omega_i }  \]
where $p_i$ and $p^\text{exact}_i$ are the pressures obtained from numerical simulation and exact solution, respectively for the $i^\text{th}$ cell whose area is $\Omega_i$. The exact pressure jump $\Delta p_\text{exact} = 36.5$ is given by the Young-Laplace law. Due to the diffused nature of the phase-field variable, the exact solution for pressure, governed by Young-Laplace law, would also be diffused. Therefore, the exact solution of pressure is defined as
\begin{equation}\label{eq:res_drop_pexact}
	p^\text{exact}_i = p_\text{min} + \left( \dfrac{\psi_i - \psi_\text{min}}{\psi_\text{max} - \psi_\text{min}} \right) \Delta p_\text{exact}
\end{equation}
For consistent comparison with the numerical solution, the minimum pressure is considered from the numerical pressure field

\[ p_\text{min} = \min_{i} p_i \]
and the same for the phase-field variable

\[ \psi_\text{min} = \min_{i} \psi_i \quad\text{and}\quad \psi_\text{max} = \max_{i} \psi_i \]

Being a static test case, the results should retain uniform zero velocity in the domain throughout the simulation. However, numerical errors in the simulation can induce spurious velocities in the domain and is a measure of the discretization errors in the numerical scheme. Along with the error in pressure field, the error in velocity field is also considered, which is measured as

\[ |\mathbf{v}|_\text{max} = \max_i |\mathbf{v}_i|\]
The two errors are tabulated in table \ref{tab:res_drop_convergence}. The $L_2$ norms of error in pressure converges with an \emph{effective} order of 1.13 for the first-order path-conservative scheme while the velocity error also decreases with grid refinement, albeit their order of magnitude remain constant.

\begin{table}[H]
	\centering
	\begin{tabular}{| c | c | c |}
		\hline
		\makecell{Mesh size (m)} & $L_2 (p)$ & $|\mathbf{v}|_\text{max}$ (m/s) \\
		\hline\hline
		$L/16$ & 0.3098 & $8.35 \times 10^{-4}$ \\
		\hline
		$L/32$ & 0.1269 & $7.64 \times 10^{-4}$ \\
		\hline
		$L/64$ & 0.0562 & $6.36 \times 10^{-4}$ \\
		\hline
		$L/128$ & 0.0300 & $5.02 \times 10^{-4}$ \\
		\hline
	\end{tabular}
	\caption{Pressure and velocity errors from static drop simulations ($L= 8$ m).}
	\label{tab:res_drop_convergence}
\end{table}

\begin{figure}[h!]
	\centering
	\includegraphics[width=\textwidth]{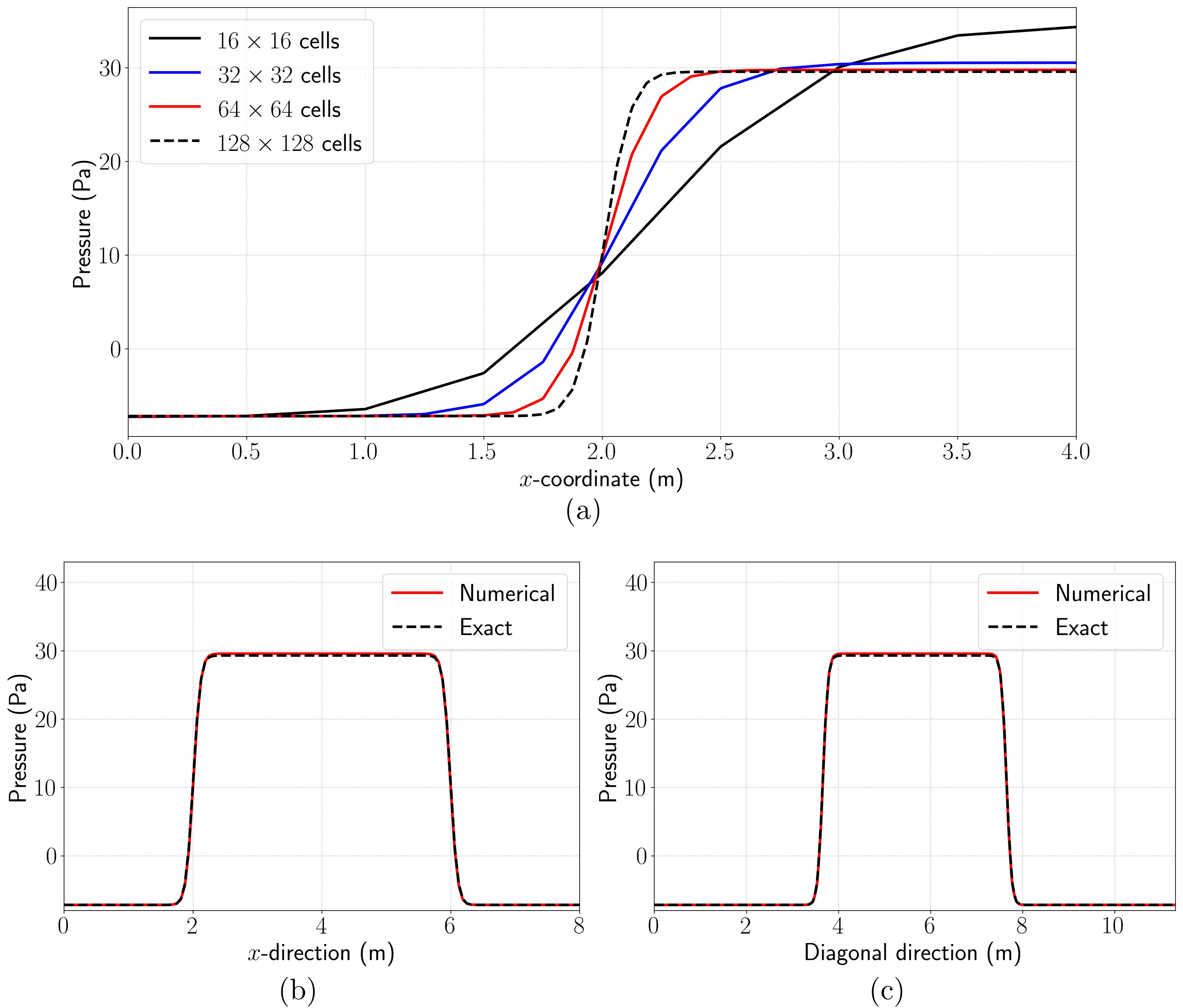}
	\caption{Static drop pressure profiles along: (a) horizontal centreline from various meshes, (b) horizontal centreline from $128 \times 128$ mesh compared against exact profile, and (c) diagonal direction from $128 \times 128$ mesh compared against exact profile.}
	\label{fig:res_drop_pressures}
\end{figure}

The pressure obtained from the various meshes, along the horizontal centreline is plotted in figure \ref{fig:res_drop_pressures}. The figure also compares the numerical pressure jump obtained from the $128 \times 128$ mesh, along the horizontal centreline as well as the diagonal direction, with the exact solution \eqref{eq:res_drop_pexact}.

\begin{figure}[ht!]
	\centering
	\includegraphics[width=4.25in]{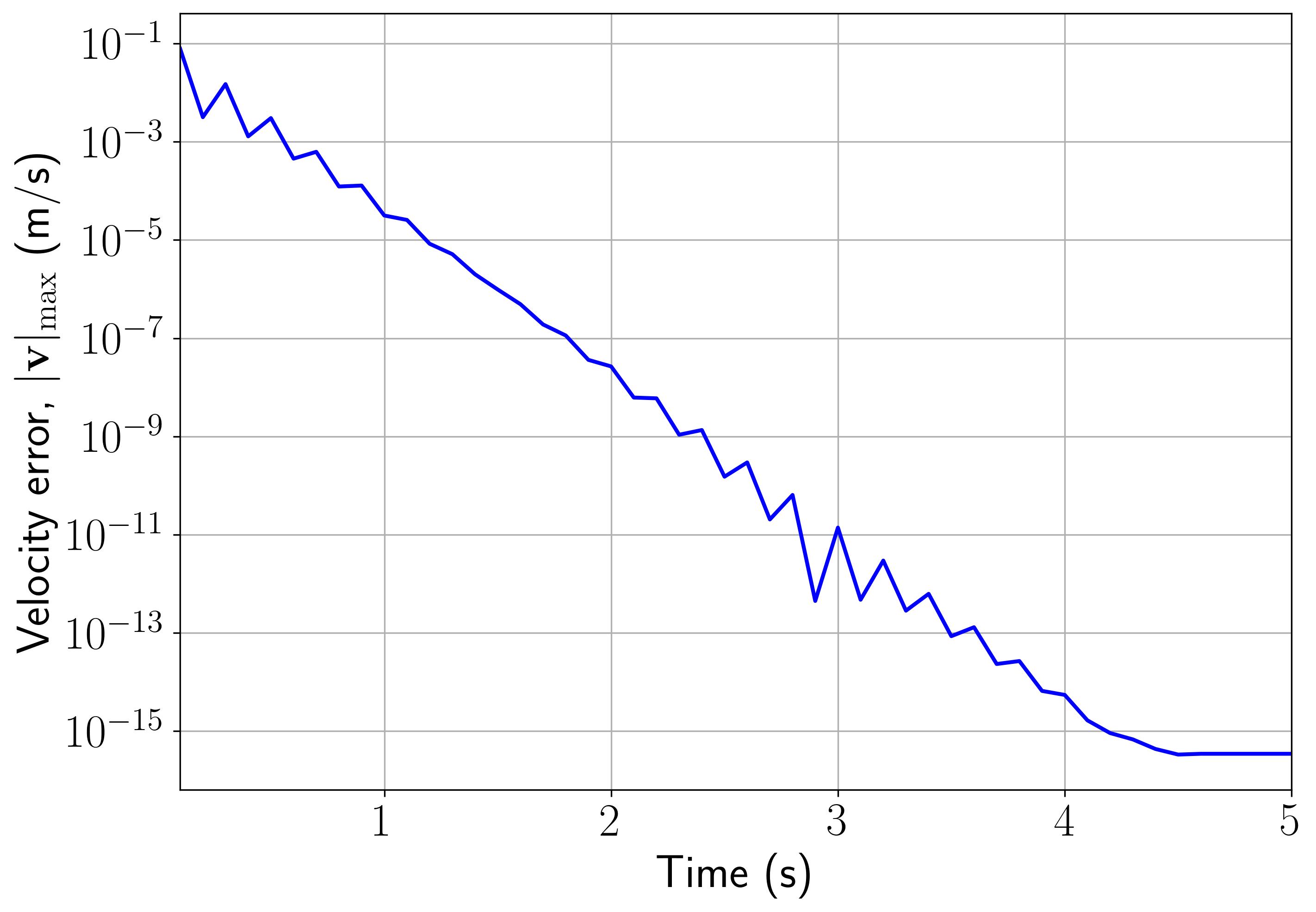}
	\caption{Temporal variation of spurious velocity in the static drop simulation on a $32 \times 32$ mesh, with the exact interface curvature prescribed.}
	\label{fig:res_drop_vmax}
\end{figure}

As previously mentioned, the presence of spurious currents in the domain is an indicator of discretization errors in the numerical scheme. However, it has been demonstrated that interface curvature computation contributes significantly to these spurious currents \cite{francois2006balanced}. Therefore, the static drop simulation is repeated on the $32 \times 32$ mesh with its exact curvature $\kappa_\text{exact} = 0.5$ prescribed. The temporal variation in velocity field error $|\mathbf{v}|_\text{max}$, plotted in figure \ref{fig:res_drop_vmax}, falls to the machine epsilon for double precision floating-point operations. This confirms that the presence of spurious currents in the previous simulations are due to the errors in interface curvature computation and not due to the discretization errors in the path-conservative scheme. This numerical exercise also verifies the well-balanced property \cite{pares2006numerical, castro2017well} of the path-conservative method developed in this work. 

To further demonstrate that the present scheme is well-balanced, along with specifying the exact curvature, the domain is initialized with exact solution of pressure \eqref{eq:res_drop_pexact} with $p_\text{min} = 0$. Simulations are conducted on a uniform $32 \times 32$ structured mesh as well as a \emph{fairly} uniform unstructured mesh with 1054 triangular cells. Throughout the simulation, the velocity field error $|\mathbf{v}|_\text{max}$ remains of the order of $10^{-17}$ m/s for structured mesh and $10^{-14}$ m/s for the unstructured mesh, showing that the path-conservative scheme is indeed well-balanced.

\subsection{Linear sloshing}

\begin{figure}[ht!]
	\centering
	\includegraphics[width=1.75in]{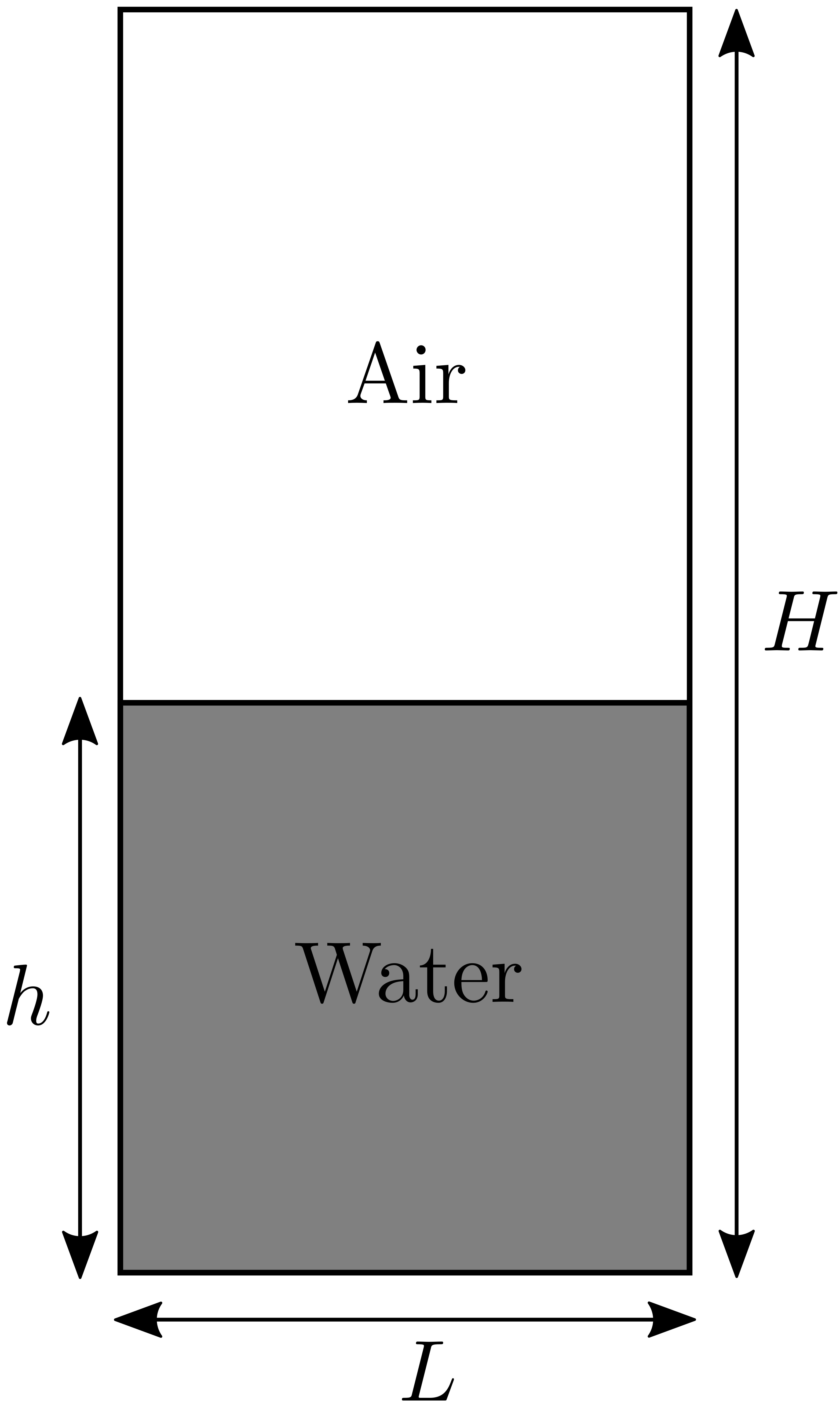}
	\caption{Schematic of the linear sloshing in a rectangular tank ($L = 1$ m, $H = 2.25$ m, and $h = 1$ m).}
	\label{fig:res_linslosh_config}
\end{figure}

Linear sloshing in a stationary rectangular tank models water sloshing under gravity and constant transverse acceleration. This problem has been studied using low Mach compressible solvers \cite{grenier2013accurate} as well as weakly compressible models \cite{li2022simplified}. The problem consists of a rectangular tank of width $L = 1$ m and height $H = 2.25$ m. The tank is filled with water till a height of $h = 1$ m and rest of the tank is filled with air. The density of water $\rho_1$ and air $\rho_2$ are $1000$ kg/m$^3$ and $1$ kg/m$^3$ respectively. The acceleration due to gravity is taken as $g_y = -9.81\,\mathrm{m/s^2}$ and the transverse acceleration is taken as $a_x = -9.81 \times 10^{-2} \,\mathrm{m/s^2}$. The initial velocity is set to zero and the pressure is set based on hydrostatic condition. Free-slip boundary condition is imposed on all four walls. The initial configuration of the problem is shown in figure \ref{fig:res_linslosh_config} and the simulations are carried out till $t = 4$ s. Being an inviscid problem, linear sloshing is an apt test case to examine the fidelity of the proposed path-conservative scheme in the absence of surface tension.

For a small acceleration ratio $a_x/g_y = 0.01$, the evolution of the  free-surface elevation $\eta$ can be computed analytically using the linearized potential theory \cite{landau1987fluid, chanteperdrix2002compressible} as

\[ \eta(x,t) = \dfrac{a_x}{g_y} \left[ \dfrac{L}{2} - x + \sum_{n \geq 0} \dfrac{4}{k^2_{2n+1} L} \cos\{\omega_{2n+1}t\} \cos\{k_{2n+1}(x-L)\} \right] \quad \forall \, x \in [0,L] \]
where

\[ k_n = \dfrac{n\pi}{L} \quad \text{and} \quad \omega_{n}^2 = \dfrac{g_y k_n (\rho_1 - \rho_2)}{\rho_1 \coth\{k_n h\} + \rho_2 \coth\{k_n (H-h)\}} \]

\begin{figure}[ht!]
	\centering
	\includegraphics[width=5.5in]{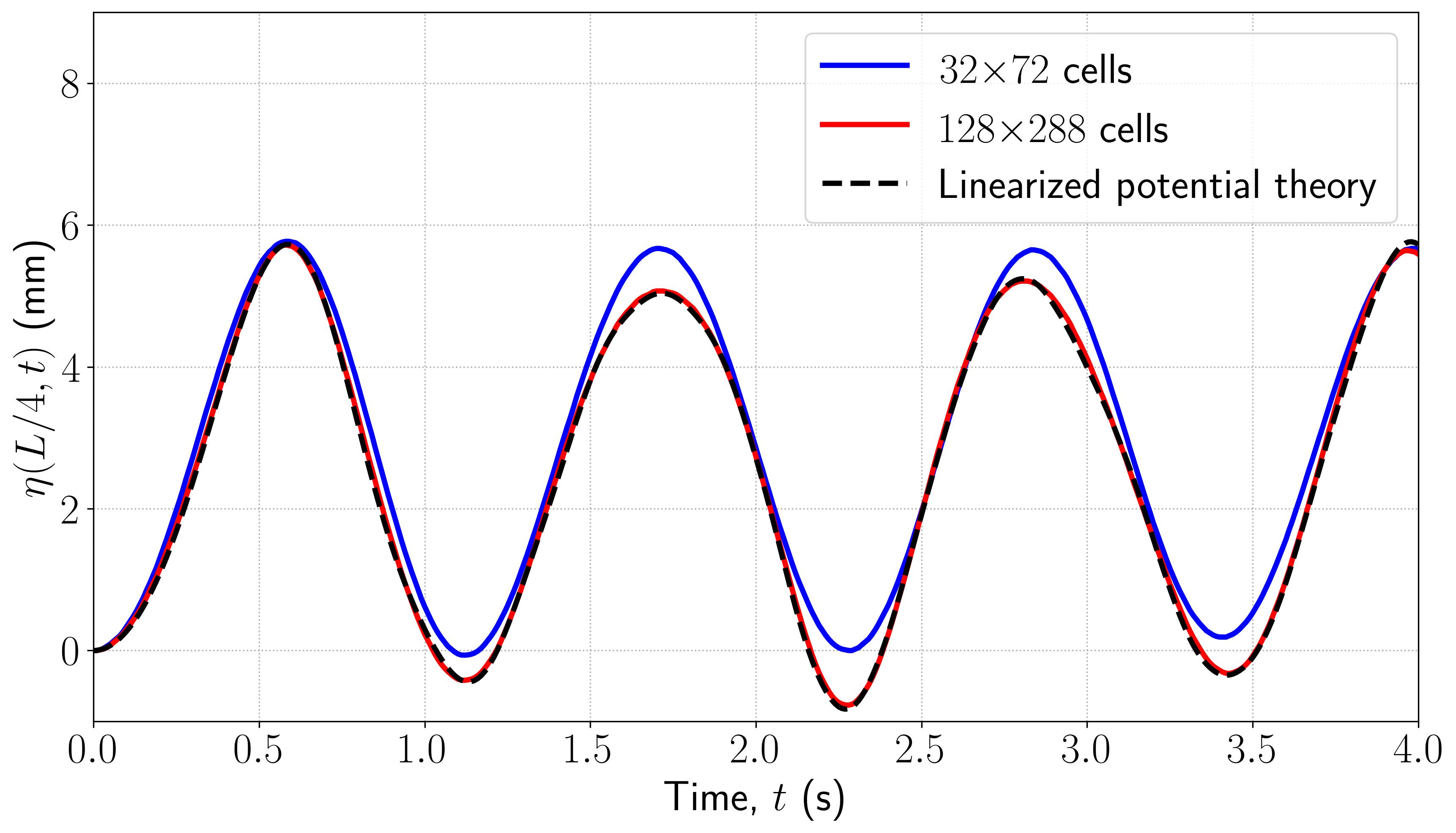}
	\caption{Elevation of air-water interface at the $x = L/4$ compared against analytical solution \cite{chanteperdrix2002compressible}.}
	\label{fig:res_linslosh_str_analytical}
\end{figure}

The numerical simulations are carried out on three uniform Cartesian meshes of sizes: $32 \times 72$, $64 \times 144$ and $128 \times 288$. The AC parameter is taken as $\beta = 1\times 10^3$. Evolution of the free-surface elevation at $x = L/4$ (a quarter of the way horizontally from the left wall), obtained from the simulations on the coarsest and the finest mesh, are compared against the analytical solution in figure \ref{fig:res_linslosh_str_analytical}. The $L_2$ norm of the errors in $\eta(L/4,t)$ with respect to the analytical solution, obtained from the three meshes are tabulated in table \ref{tab:res_linslosh_convergence}. The path-conservative scheme with linear reconstruction yields an \emph{effective} convergence order of 1.39 for the linear sloshing problem. 

\begin{table}[H]
	\centering
	\begin{tabular}{| c | c | c |}
		\hline
		\makecell{Mesh size (m)} & $L_2$ error & Order \\
		\hline\hline
		1/32 & 0.00825 & - \\
		\hline
		1/64 & 0.00374 & 1.14 \\
		\hline
		1/128 & 0.00119 & 1.65 \\
		\hline
	\end{tabular}
	\caption{Convergence order of error in linear sloshing.}
	\label{tab:res_linslosh_convergence}
\end{table}

\begin{figure}[ht!]
	\centering
	\includegraphics[width=4.5in]{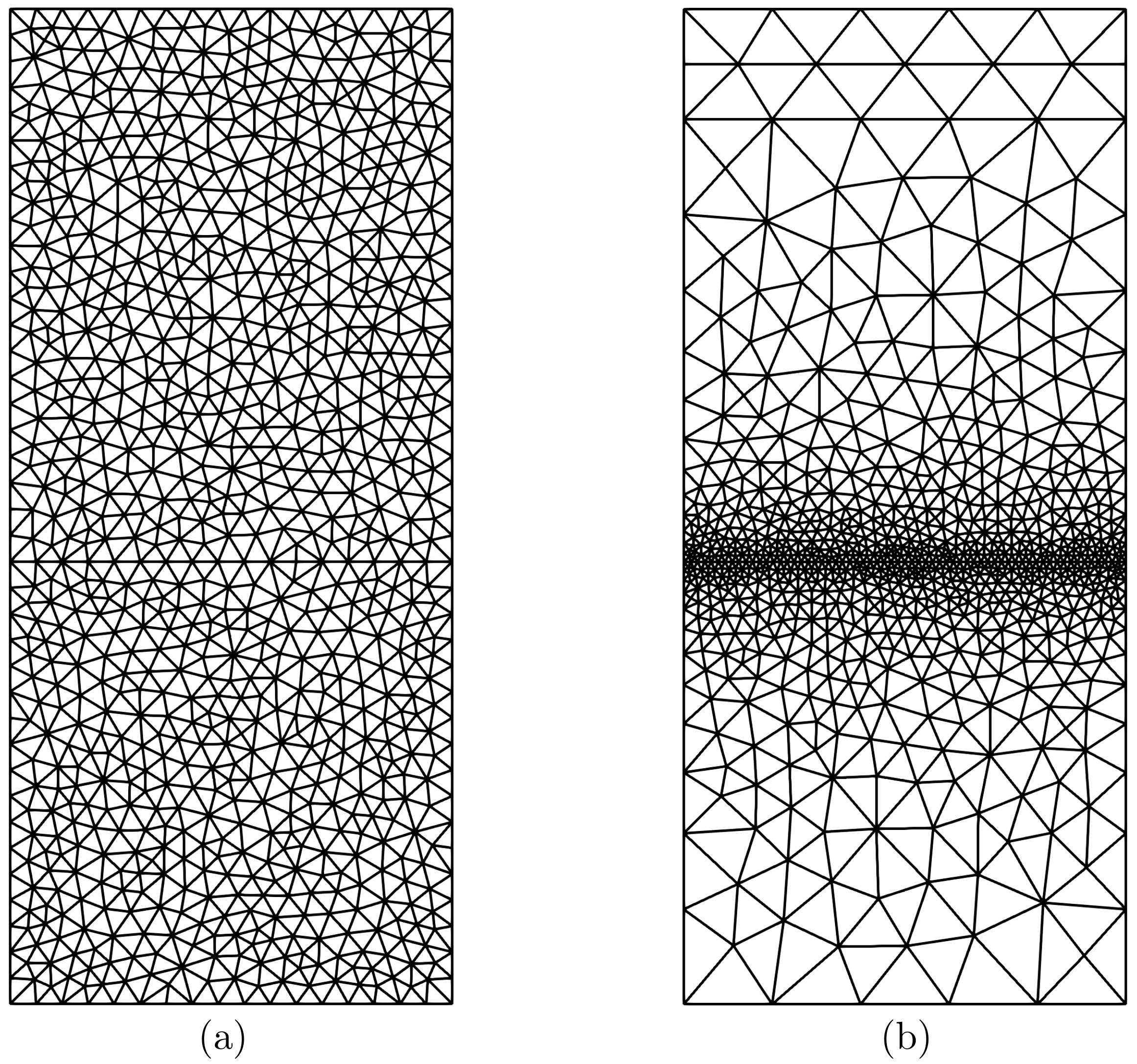}
	\caption{Unstructured meshes used for linear sloshing problem (approx. 1700 triangular cells): (a) uniform mesh, and (b) mesh refined near the interface.}
	\label{fig:res_linslosh_unstr_meshes}
\end{figure}

\begin{figure}[ht!]
	\centering
	\includegraphics[width=5.5in]{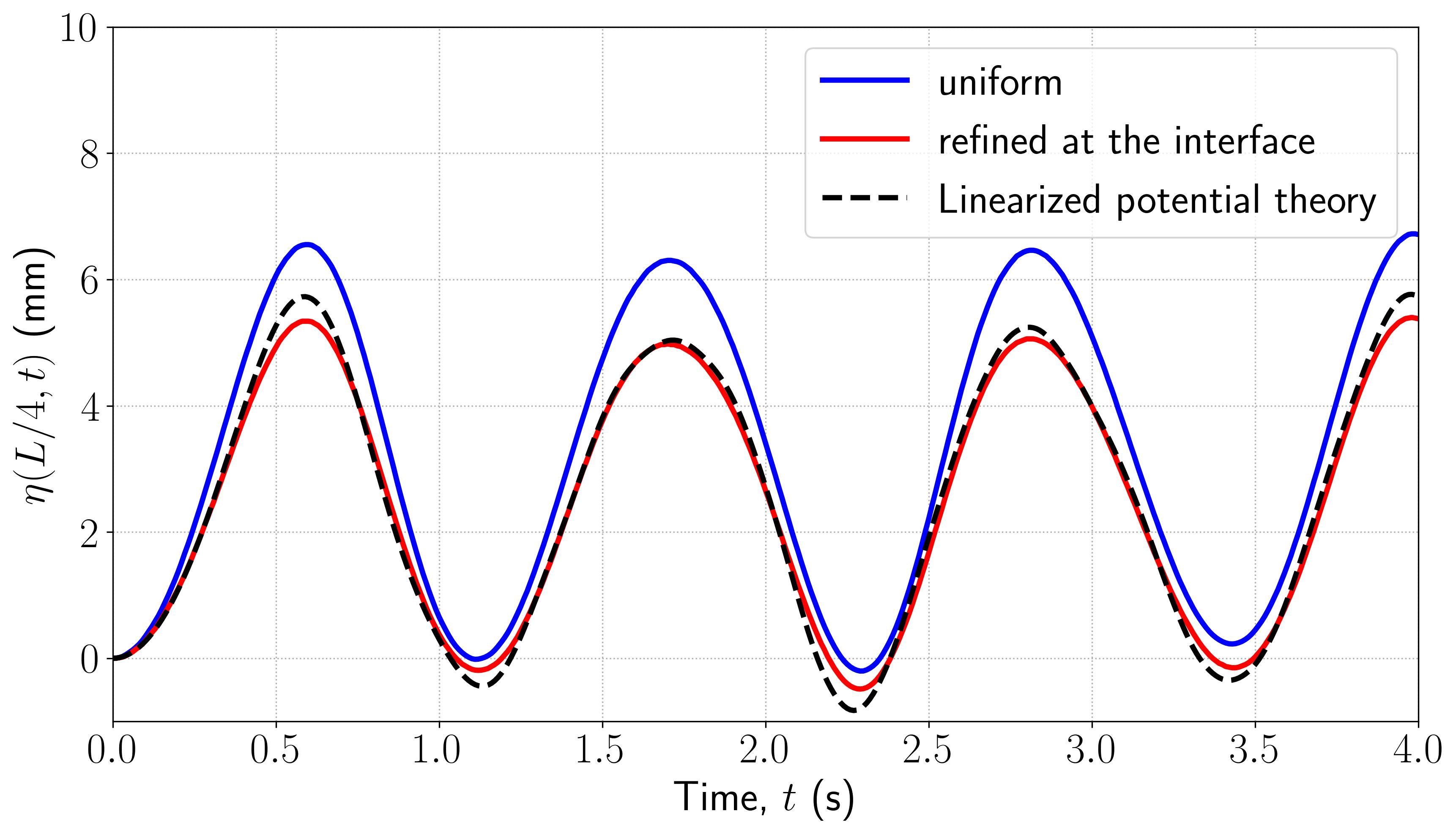}
	\caption{Elevation of air-water interface at $x = L/4$, obtained from unstructured meshes, compared against analytical solution \cite{chanteperdrix2002compressible}.}
	\label{fig:res_linslosh_unstr_analytical}
\end{figure}

To demonstrate the adaptability of the developed solver, linear sloshing problem is solved on unstructured meshes with triangular cells. Two different meshes with approximately the same number of cells were used: a) (fairly) uniform sized triangular cells across the entire domain, and b) triangular cells refined near the region of interest (interface), as shown in figure \ref{fig:res_linslosh_unstr_meshes}. The evolution of the fluid interface at $x = L/4$, obtained from the unstructured mesh simulations, is compared against analytical solution in figure \ref{fig:res_linslosh_unstr_analytical}, demonstrating the capability of the solver to handle arbitrary unstructured meshes. Predictably, the mesh refinement near the interface yields more accurate results with its $L_2$ error norm comparable to the one from $64 \times 144$ uniform Cartesian mesh. The ability of the developed solver to handle non-uniform meshes aids in accurate interface capturing with reduced computational overhead, which can be particularly desirable in large-scale problems.

\subsection{Capillary wave}

The damping of a capillary wave is an ideal validation test case as it requires modelling viscous fluxes along with the inclusion of surface tension effects in the hyperbolic solver. Furthermore, the oscillation of the damped wave has an analytical solution \cite{cortelezzi1981small, prosperetti1981motion}. The amplitude of the capillary wave, over time, is given by the exact solution to the initial-value problem of small-amplitude waves between two viscous fluids \cite{prosperetti1981motion}.

\begin{figure}[ht!]
	\centering
	\includegraphics[width=3in]{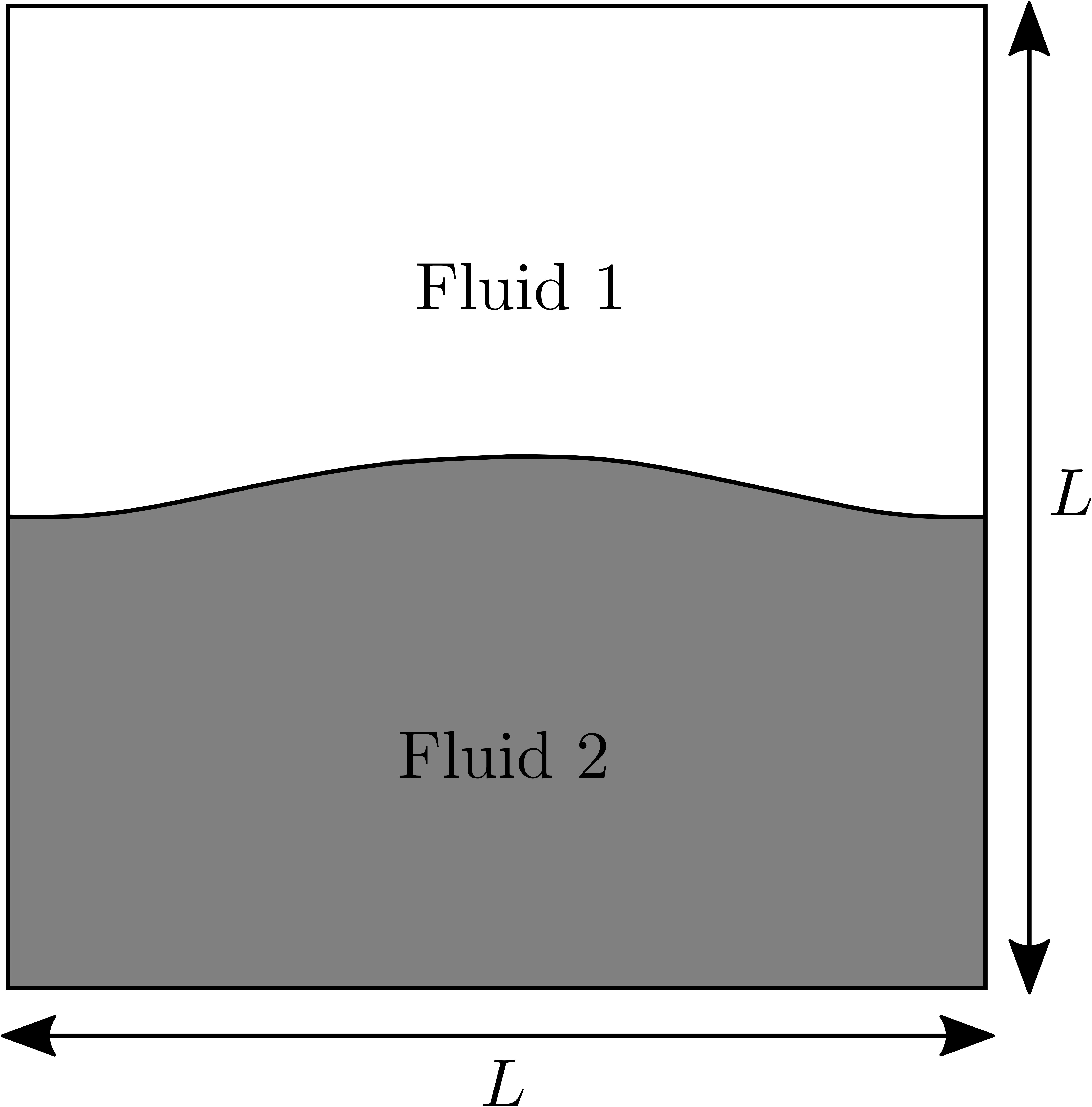}
	\caption{Schematic of the capillary wave problem ($L = 1$ m).}
	\label{fig:res_capwave_config}
\end{figure}

A $L \times L$ square domain, where $L = 1$ m is considered for this problem. Two fluids of the same density $\rho_1 = \rho_2 = 18.3\,\mathrm{kg/m}^3$ and dynamic viscosity $\mu_1 = \mu_2 = 7.8 \times 10^{-2} \,$kg/m-s are initially separated by an interface defined as 

\[ y(x) = L/2 - A_0 \cos (2\pi x) \quad \forall \, x \in [0,L]  \] 
where the amplitude $A_0 = 0.01$ m. The surface tension between the two fluids is taken as $\sigma = 1$ N/m. Initial configuration of the problem is shown in figure \ref{fig:res_capwave_config}. The left and right boundaries form a pair of periodic boundaries while the top and bottom boundaries are considered as free-slip walls. The initial pressure and velocity are set to zero.

\begin{figure}[ht!]
	\centering
	\includegraphics[width=5.5in]{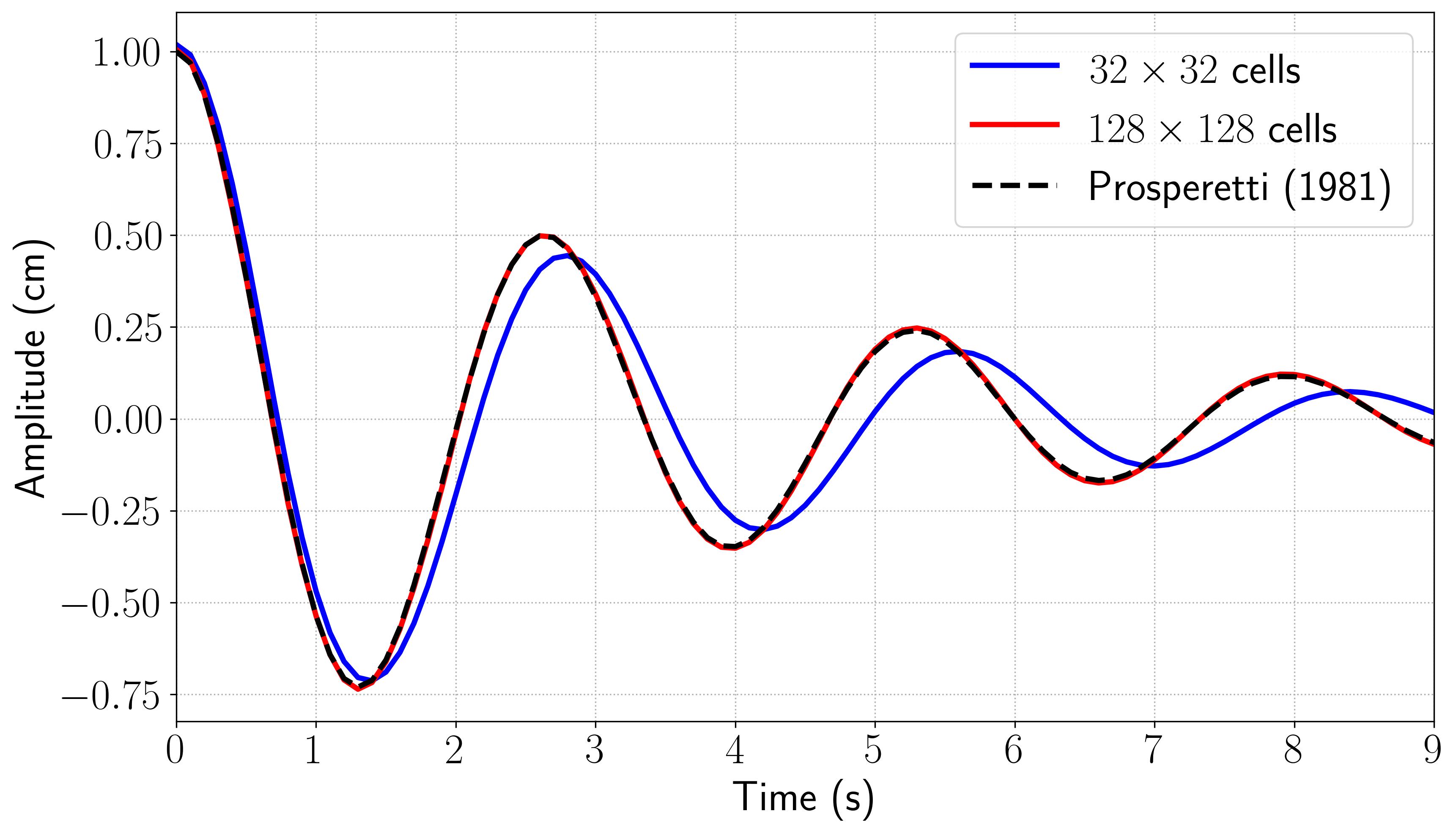}
	\caption{Amplitude of the capillary wave, measured at $x = L/2$, compared against analytical solution \cite{prosperetti1981motion}.}
	\label{fig:res_capwave_analytical}
\end{figure}

The simulations are carried out on three uniform Cartesian meshes of sizes: $32 \times 32$, $64 \times 64$, and $128 \times 128$, till $t = 9$ s, with AC parameter $\beta = 50$. The temporal variation of amplitude of the wave obtained from the coarsest and the finest mesh simulations are compared against the analytical solution in figure \ref{fig:res_capwave_analytical}. The $L_2$ norm of the errors in amplitude with respect to the analytical solution, obtained from the three meshes are tabulated in table \ref{tab:res_capwave_convergence}.  The path-conservative scheme with linear reconstruction yields an \emph{effective} convergence order of 2.11 for this case.

\begin{table}[H]
	\centering
	\begin{tabular}{| c | c | c |}
		\hline
		\makecell{Mesh size (m)} & $L_2$ error & Order \\
		\hline\hline
		1/32 & 0.00962 & - \\
		\hline
		1/64 & 0.00256 & 1.91 \\
		\hline
		1/128 & 0.00051 & 2.32 \\
		\hline
	\end{tabular}
	\caption{Convergence order of error in capillary wave.}
	\label{tab:res_capwave_convergence}
\end{table}

To demonstrate the need to include surface tension in the hyperbolic solver, the present solver is compared against an alternate HLLC path-conservative scheme. In this alternate solver, the surface tension terms are ignored while formulating the path-conservative scheme in section \ref{sec:num}. The surface tension terms are modelled separately using the continuum surface force (CSF) model \cite{brackbill1992continuum} and treated as a source term \cite{parameswaran2023stable}. The formulation of this alternate (HLLC+CSF) solver, follows the same procedure as in section \ref{sec:num}. The key difference in HLLC+CSF formulation is that the generalized Riemann invariants across contact wave shows no pressure jump across it, unlike in the present work \eqref{eq:num_laplacelaw}. Therefore, the intermediate pressures are equal ($p_{*L} = p_{*R}$) in the HLLC formulation \eqref{eq:num_hllcStar} of HLLC+CSF solver.

\begin{figure}[ht!]
	\centering
	\includegraphics[width=5.5in]{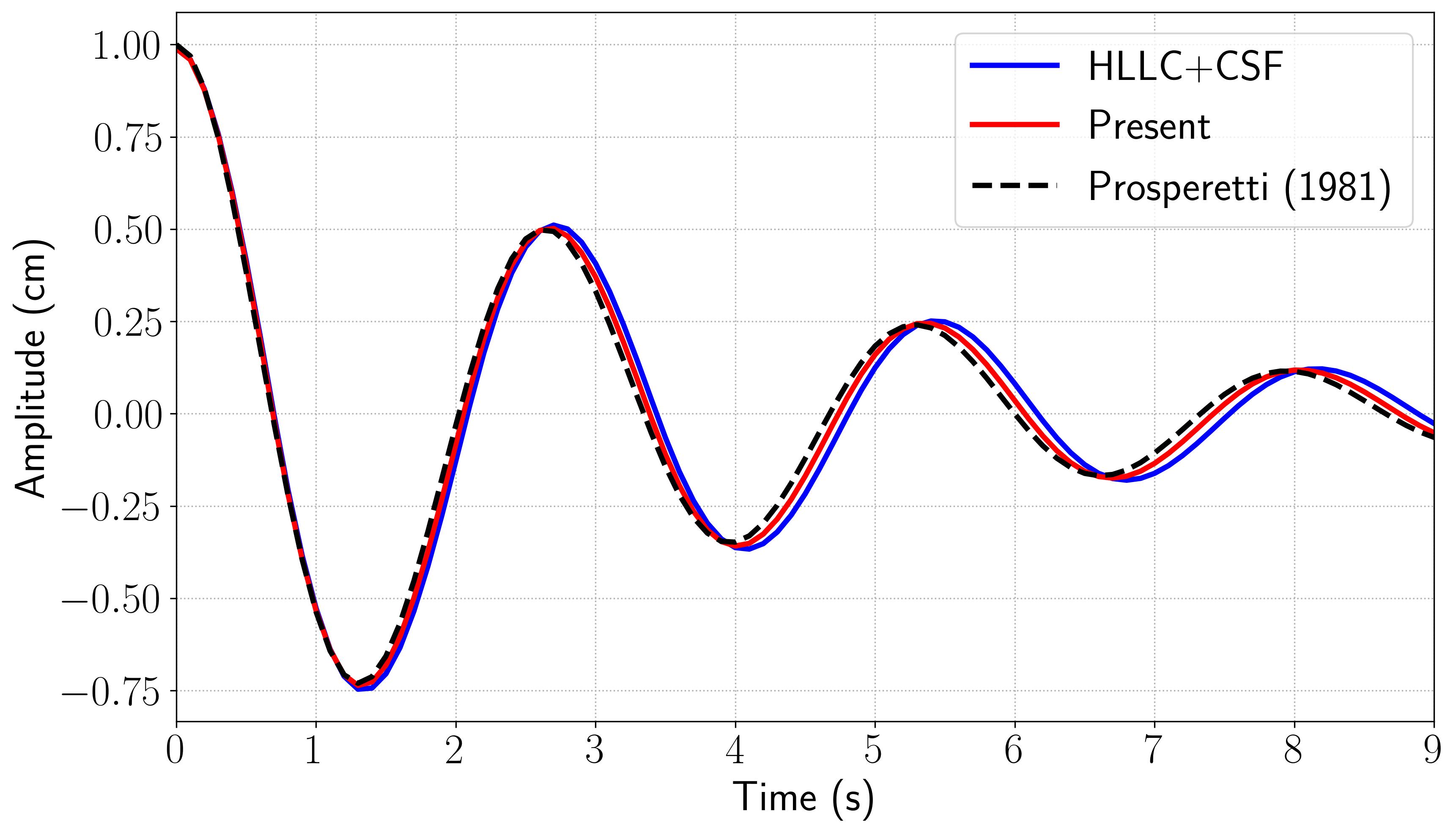}
	\caption{Amplitude of the capillary wave obtained from the present solver compared against an alternate hyperbolic solver (HLLC+CSF) without surface tension terms in the Riemann solver.}
	\label{fig:res_capwave_hllc_csf}
\end{figure}

\begin{figure}[ht!]
	\centering
	\includegraphics[width=0.7\textwidth]{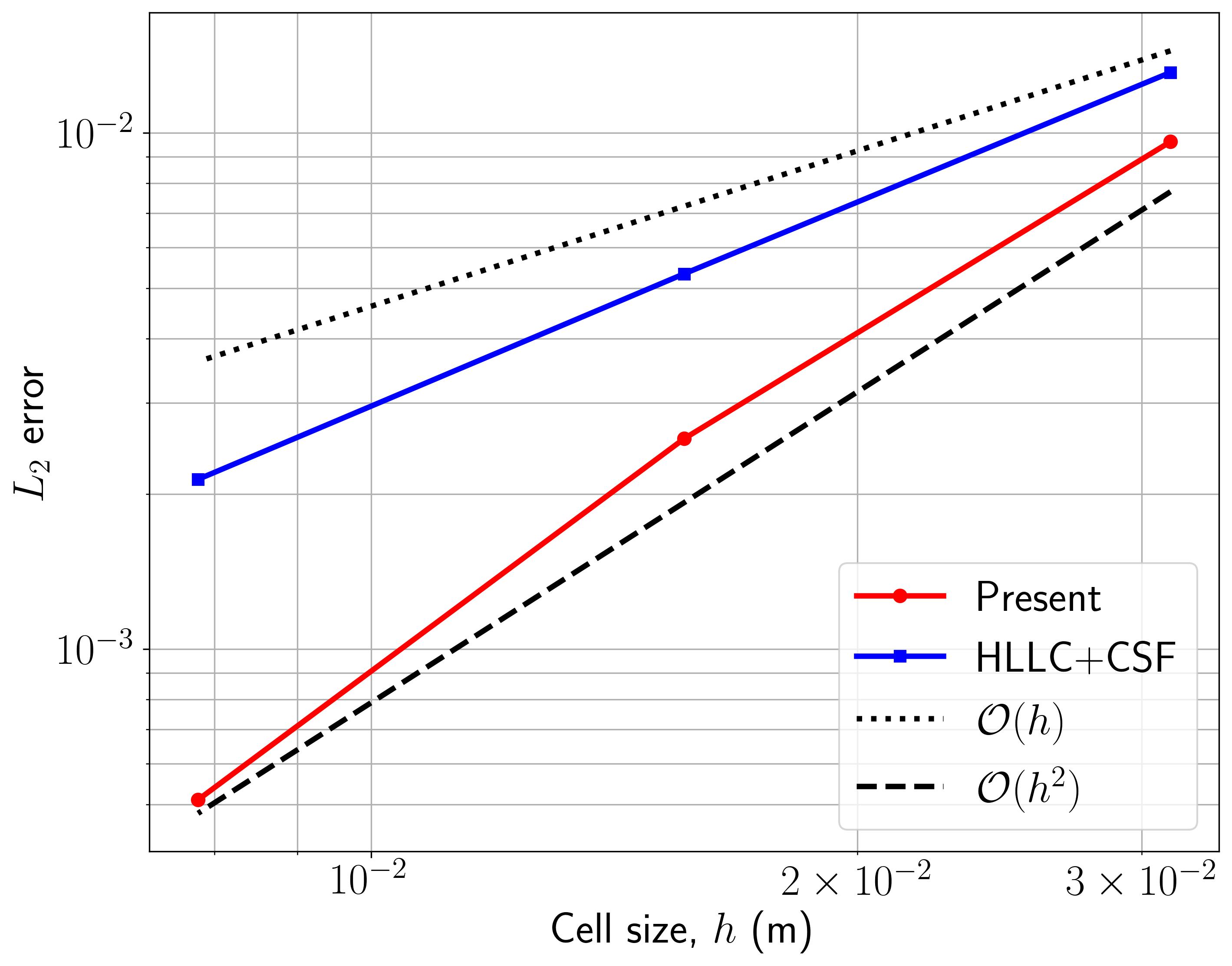}
	\caption{Convergence of $L_2$ norm of error in the amplitude of the capillary wave obtained from the present solver compared against an alternate hyperbolic solver (HLLC+CSF) without surface tension terms in the Riemann solver.}
	\label{fig:res_capwave_hllc_csf_convergence}
\end{figure}

The capillary wave problem is simulated using the two approaches on a $64 \times 64$ uniform Cartesian mesh and the temporal evolution of the wave amplitude is compared against analytical solution in figure \ref{fig:res_capwave_hllc_csf}. The results from HLLC+CSF solver deviates from the analytical solution which becomes more pronounced over time. The $L_2$ norm of the errors from HLLC+CSF solver is compared against that of the present solver in figure \ref{fig:res_capwave_hllc_csf_convergence}. The \emph{effective} convergence order obtained using HLLC+CSF solver is 1.31, which is considerably lower than that obtained from the present solver. The numerical exercise clearly indicates that decoupling the surface tension terms from the hyperbolic solver can lead to inaccuracies in surface tension dominated flows where the present framework can improve the results significantly.

\subsection{Rayleigh-Taylor instability}

The Rayleigh-Taylor instability (RTI) considers the interface between two fluids when a heavier fluid rests over a lighter fluid and gravity is acting in the direction from heavier to lighter fluid. The instability grows exponentially in time, if the interface is disturbed by a small perturbation. For inviscid incompressible fluids, the growth rate $n$ of the instability is given by the linear theory \cite{bellman1954effects}

\begin{equation}\label{eq:res_rti_n}
	n^2 = k\, |g_y| \left[ A - \dfrac{k^2 \sigma}{|g_y|(\rho_1 + \rho_2)} \right]
\end{equation}
Here $k$ is the wavenumber of the initial perturbation, $g_y$ is the gravitational acceleration perpendicular to the interface and $\sigma$ is the surface tension coefficient. Non-dimensional $A$ is the Atwood number defined as $A = (\rho_1 - \rho_2)/(\rho_1 + \rho_2)$, where $\rho_1$ and $\rho_2$ are the densities of the heavier and lighter fluids respectively. Following \cite{daly1967numerical}, to study the effect of surface tension on the growth rate, it would be helpful to introduce the parameter $\phi_s$ defined as

\[ \phi_s = \dfrac{k^2 \sigma}{A|g_y|(\rho_1 + \rho_2)} \]

\begin{figure}[h!]
	\centering
	\includegraphics[width=1.75in]{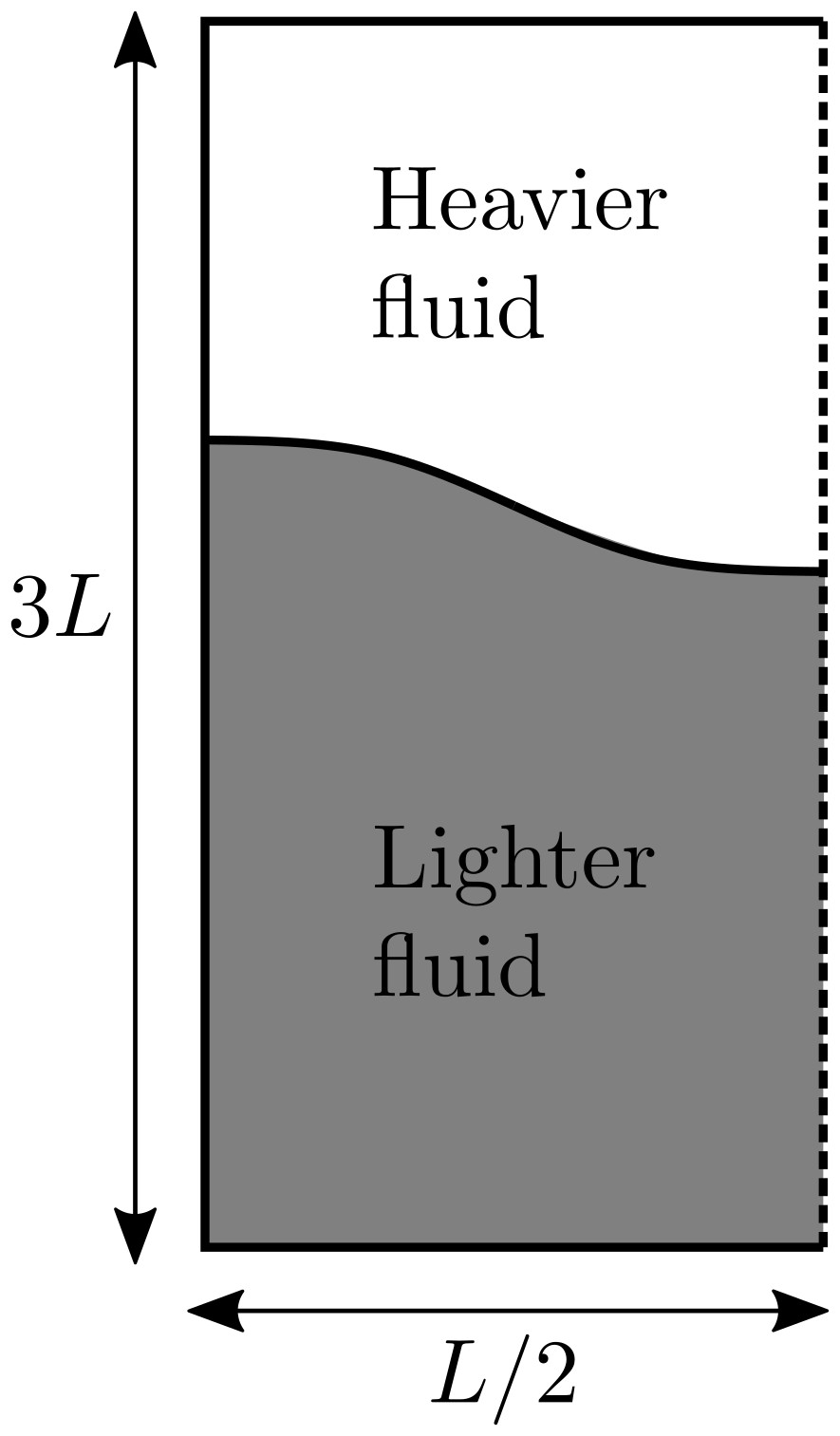}
	\caption{Schematic of the Rayleigh-Taylor instability problem (figure not to proportion). The right dashed wall indicates the line of symmetry ($L = 6 \pi /5$ m).}
	\label{fig:res_rti_schematic}
\end{figure}

For the test case, the densities of the heavier and lighter fluids are taken as $\rho_1 = 1.2\,\mathrm{kg/m}^3$ and $\rho_2 = 0.3\,\mathrm{kg/m}^3$ respectively resulting in Atwood number $A = 0.6$. The gravitational acceleration is taken as ${\mathbf{g}} = (0,-1)\,\mathrm{m/s^2}$. The wavenumber of the initial perturbation is taken as $k = 5/3$ leading to a column width $L = 2\pi/ k  = 6\pi/5$. Owing to the symmetry of the problem, only half-width is considered leading to a $L/2 \times 3L$ computational domain. The initial interface location is prescribed by

\[ y(x) = 2L + 0.035\cos (k x) \quad \forall \, x \in [0,L/2]  \]
with $x = L/2$ being the line of symmetry. All the other boundaries are taken as free-slip walls. The initial configuration of the problem is shown in figure \ref{fig:res_rti_schematic}.

\begin{figure}[ht!]
	\centering
	\includegraphics[width=0.9\textwidth]{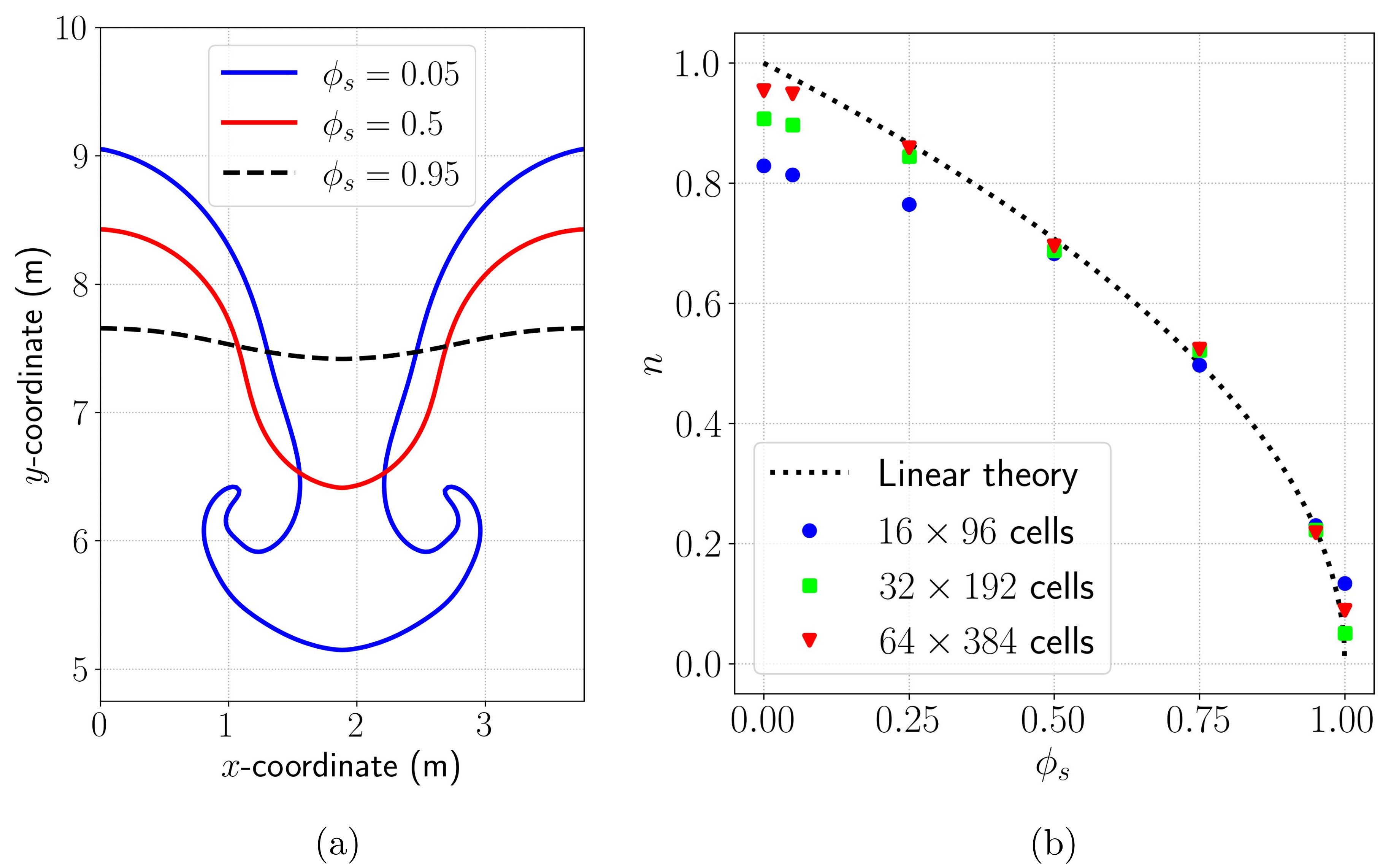}
	\caption{Results from Rayleigh-Taylor instability simulations: (a) Interface shapes obtained from $64 \times 384$ mesh at $t = 6$ s for three different values of $\phi_s$, and (b) Growth rate obtained from three different meshes for various $\phi_s$ compared against linear theory \cite{bellman1954effects}.}
	\label{fig:res_shape_nPhi}
\end{figure}

The RTI problem is solved in three uniform Cartesian meshes of sizes: $16 \times 96$, $32 \times 192$ and $64 \times 384$. The simulations are carried out till $t = 10$ s and the AC parameter is taken as $\beta = 1 \times 10^3$. To study the effect of surface tension on the instability, parameter $\phi_s$ is varied between 0 and 1. The interface shapes at $t = 6$ s obtained from the fine mesh simulation, for three different $\phi_s$, are shown in figure \ref{fig:res_shape_nPhi}. In the figure, interface profiles along the entire width $L$ is plotted by mirroring the half-width profile about the line of symmetry. It is evident from the figure that the growth of instability is delayed as surface tension increases.

To estimate the growth rate $n$, the temporal variation of logarithm of amplitude, at the line of symmetry, is plotted. In the plot, a linear region is identified and its slope is estimated by linear regression. The slope gives the estimated growth rate $n$ of the instability. This process is repeated for various values of $\phi_s \in [0,1]$ and compared against analytical solution \cite{bellman1954effects, daly1967numerical} in figure \ref{fig:res_shape_nPhi}. The present solver can capture the growth rate accurately in the fine mesh simulations for all values $\phi_s$. As $\phi_s$ increases, the growth of the instability is delayed and the results from all three meshes coincide. It is worth pointing out the deviation in the results at $\phi_s = 1$, where the surface tension stabilizes the interface leading to no growth in instability ($n = 0$). The deviation may be attributed to the spurious currents induced by the computation of interface curvature, as discussed in the case of static drop. The deviation may also be the result of unbalanced forces across the fluid interface due to gravity, which needs to be further verified.

The growth rates obtained are comparable to the ones obtained from incompressible solvers reported in \cite{gerlach2006comparison}. On the other hand, the growth rates observed in incompressible two-phase simulations of Rayleigh-Taylor instability using artificial compressibility method \cite{bhat2022improved} deviate significantly from linear theory when surface tension effects are considered. In \cite{bhat2022improved}, the surface tension is modelled using the continuum surface stress (CSS) \cite{lafaurie1994modelling} method, independent of the hyperbolic solver and therefore overlooks the jump in pressure across the interface, further substantiating the importance of including surface tension effects in the hyperbolic system.

\section{Conclusion}\label{sec:conc}

In this work, a HLLC-type path-conservative method has been developed for a weakly compressible (WC) two-phase model. The proposed scheme offers two advantages over our previous method \cite{melvin2025development}: a) the non-conservative terms are treated by the path-conservative scheme, eliminating the need for a separate discretization technique, and b) the surface tension effects are included in the non-conservative system resulting in accurate jump conditions and improved accuracy, leading to a more robust formulation. The first-order path-conservative scheme is combined with a local solution reconstruction technique to obtain high-order spatial accuracy.

The efficacy of the HLLC-type path-conservative scheme is tested on several incompressible two-phase flow problems. The case of static drop was used to demonstrate the capability of the path-conservative scheme to accurately capture the pressure jump across a circular fluid interface due to surface tension. The results also confirmed the well-balanced property of the present scheme. The simulations of inviscid linear sloshing demonstrated the fidelity of the proposed path-conservative method in the absence of surface tension. Since the proposed finite volume formulation is for an arbitrary two-dimensional mesh, the linear sloshing problem was tested on unstructured meshes as well. The case of viscous damping of a capillary wave was utilized to highlight the need to include surface tension terms in the non-conservative system of the WC two-phase model. Finally, the effect of surface tension on the growth of the Rayleigh-Taylor instability was studied using the proposed solver with the obtained growth rates agreeing well with linear theory. 

The spatial orders of convergence were also reported for the first-order as well as high-order schemes, which were on par with the expected nominal orders for the respective schemes. Spatial order obtained without the inclusion of surface tension terms in the hyperbolic system was found to be considerably lower which further highlights the contentious treatment of surface tension terms as a source term. Inclusion of other source terms, such as gravitational force, in the hyperbolic system needs to be studied which can enhance the accuracy of WC two-phase models.


\newpage


\printbibliography[title=References]

\end{document}